\documentclass[12pt]{article}

\usepackage{graphicx}
\usepackage{grffile}
\usepackage{amsmath}
\usepackage{amssymb}
\usepackage{amstext}
\usepackage{amscd}
\usepackage{amsfonts}
\usepackage{xcolor}
\usepackage{wasysym}
\usepackage{stackrel}
\usepackage{stackengine}
\usepackage{txfonts}

\usepackage{a4wide}
\usepackage{authblk}

\providecommand{\keywords}[1]{\textbf{Keywords: } #1}

\newcommand\owedge{\stackMath\mathbin{\stackinset{c}{0ex}{c}{0ex}{\scriptstyle{\wedge}}{\medcirc}}}

\numberwithin{equation}{section}

\begin{document}


\title{Along the lines of nonadditive entropies: 
       $q$-prime numbers and $q$-zeta functions}

\author[1,5]{Ernesto P.\ Borges 
             \thanks{ernesto@ufba.br} 
            }
\author[2,3]{Takeshi Kodama 
             \thanks{tkodama@if.ufrj.br}
}
\author[4,5,6,7]{Constantino Tsallis 
                 \thanks{tsallis@cbpf.br}
                }

\affil[1]{Instituto de F\'isica, Universidade Federal da Bahia,
          Rua Bar\~ao de Jeremoabo s/n, Salvador-BA 40170-115, Brazil}
\affil[2]{Instituto de F\'isica, Universidade Federal do Rio de Janeiro, 
          Rio de Janeiro-RJ  21941-972, Brazil}
\affil[3]{Instituto de F\'isica, Universidade Federal Fluminense 
          and 
          National Institute of Science and Technology 
          for Nuclear Physics and Applications, 
          Campus da Praia Vermelha, Niter\'oi-RJ  24210-346, Brazil}
\affil[4]{Centro Brasileiro de Pesquisas F\'{\i}sicas, 
          Rua Xavier Sigaud 150, 22290-180 - Rio de Janeiro, RJ, Brazil}
\affil[5]{National Institute of Science and Technology of Complex Systems, 
\newline
          Rua Xavier Sigaud 150, 22290-180 - Rio de Janeiro, RJ, Brazil}
\affil[6]{Santa Fe Institute -- 1399 Hyde Park Road, Santa Fe, 87501 NM, United States}
\affil[7]{Complexity Science Hub Vienna - Josefst\"adter Strasse 39, 
          1080 Vienna, Austria}

\date{\empty}

\maketitle

\abstract{
The rich history of prime numbers includes great names such as  Euclid, 
who first analytically studied the prime numbers 
and proved that there is an infinite number of them, 
Euler, who introduced the function 
$\zeta(s) \equiv \sum_{n=1}^\infty n^{-s} = 
\prod_{p\,prime} \frac{1}{1- p^{-s}}$, 
Gauss, who estimated the rate at which prime numbers increase,
and Riemann, who extended $\zeta(s)$ to the complex plane $z$ 
and conjectured that all nontrivial zeros are in the $\mathbb{R}(z)=1/2$ axis. 
The nonadditive entropy 
$S_q= k \sum_ip_i\ln_q(1/p_i) \;(q \in \mathbb{R}; 
\,S_{1}=S_{BG}\equiv - k\sum_i p_i \ln p_i$, 
where BG stands for Boltzmann-Gibbs) 
on which nonextensive statistical mechanics is based, 
involves the function 
$ \ln_q z \equiv \frac{z^{1-q}-1}{1-q}\; (\ln_1 z=\ln z)$. 
It is already known that
this function paves the way for the emergence of a $q$-generalized algebra, 
using $q$-numbers defined as $\langle x\rangle_q \equiv e^{\ln_q x}$, 
which recover the  number  $x$ for $q=1$. 
The $q$-prime numbers are then defined as the $q$-natural numbers
$\langle n\rangle_q \equiv e^{\ln_q n} \;(n=1,2,3,\dots)$, 
where $n$ is a prime number $p=2,3,5,7,\dots$ 
We show that, for any value of $q$, infinitely many 
$q$-prime numbers exist; 
for $q\le 1$ they diverge for increasing prime number, 
whereas they converge for $q>1$; 
the standard prime numbers are recovered for $q=1$. 
For $q\le 1$, we generalize the  $\zeta(s)$ function as follows: 
$\zeta_q(s) \equiv  \langle\zeta(s)\rangle_q$ ($s \in \mathbb{R}$). 
We show that this function appears to diverge at $s=1 +0$, $\forall q$.
Also, we alternatively define, for $q\le 1$,
$\zeta_q^{\Sigma}(s) 
  \equiv  \sum_{n=1}^\infty \frac{1}{\langle n\rangle_q^s} 
  =       1+\frac{1}{\langle 2\rangle_q^s}+ \dots $ 
and  
$\zeta_q^{\Pi}(s) 
 \equiv \prod_{p\,prime} \frac{1}{1-\langle p\rangle_q^{-s}} 
  =     \frac{1}{1-\langle 2\rangle_q^{-s}} 
        \frac{1}{1-\langle 3\rangle_q^{-s}} 
        \frac{1}{1-\langle 5\rangle_q^{-s}} \cdots$, 
which, for $q<1$, generically satisfy $\zeta_q^\Sigma(s) < \zeta_q^\Pi(s)$, 
in variance with the $q=1$ case, where of course 
$\zeta_1^\Sigma(s) = \zeta_1^\Pi(s)$.
}

\medskip

\section{Introduction}
\label{sectionintroduction}

Extending the realm of the Boltzmann-Gibbs-von Neumann-Shannon 
entropic functional, 
many measures of uncertainty have been proposed to handle complex systems 
and, ultimately, complexity.  
All of them are nonadditive, excepting the Renyi functional. 
Among them, a paradigmatic one is the entropy $S_q$ defined, 
with the scope of generalizing Boltzmann-Gibbs (BG) statistical mechanics,
as follows \cite{Tsallis1988}:
\begin{equation}
S_q=k\frac{1-\sum_i p_i^q}{q-1}=k\sum_i p_i \ln_q \frac{1}{p_i}\;\;(k>0)
\end{equation}
with $S_1=S_{BG}\equiv k\sum_ip_i\ln\frac{1}{p_i}$. 
The entropy $S_q$ is the most general one which simultaneously 
is composable and trace-form \cite{EncisoTempesta2017}, 
and it has been shown to be connected to the Euler-Riemann function 
$\zeta(s)$ \cite{Tempesta2011}. 
The $q$-logarithm function is defined \cite{Tsallis1994} as follows
\begin{equation}
 \label{eq:q-log}
 \ln_q z \equiv \frac{z^{1-q}-1}{1-q}\;\;(\ln_1 z=\ln z)\,,
\end{equation}
its inverse function being
\begin{equation}
 \label{eq:q-exp}
 e_q^z=[1+(1-q)z]^{1/(1-q)} \;(e_1^z=e^z)\,
\end{equation}
when $[1+(1-q)z] \ge 0$; otherwise it vanishes.
The definitions of the $q$-logarithm and $q$-exponential functions
allow consistent generalizations of algebras 
\cite{lemans-2003,Borges2004,lemans-2009,BorgesCosta2021,Gomez-Borges-2021},
calculus
\cite{Borges2004,Kalogeropoulos-2005,
         Nobre-RegoMonteiro-Tsallis-2011,
         BorgesCosta2021}
(see also \cite{Czachor-2020})
and generalized numbers \cite{Lobao-etal2009,BorgesCosta2021}.

There are different ways of defining generalized $q$-numbers connected with
the pair of inverse ($q$-logarithm, $q$-exponential) functions, namely
\begin{eqnarray}
\label{eq:e_lnq_x}
\langle x \rangle_q &=& e^{\ln_q x}\;\;\, (\langle x\rangle_1=x>0)\,,
\\
\label{eq:eq_ln_x}
{}_q\langle x \rangle &=& e_q^{\ln x}\;\;\;\; ({}_1\langle x \rangle =x>0)\,,
\\
\label{eq:ln_eq_x}
[x]_q &=& \ln e_q^x \;\;\; ([x]_1=x \in \mathbb{R})\,,
\\
\label{eq:lnq_e_x}
{}_q[x] &=& \ln_q e^x\;\, (_1[x]=x \in \mathbb{R})\,.
\end{eqnarray}
Observe that
$\langle\,{}_q\langle x \rangle\,\,\rangle_q = 
 {}_q\langle\:\langle x\rangle_q\,\rangle = 
 [\,{}_q[x]\,\,]_q = 
 {}_q[\:[x]_q\,] = x$.
These four possibilities are explored in Ref.\ \cite{BorgesCosta2021}%
\footnote{Eqs.\  (\ref{eq:e_lnq_x}), (\ref{eq:eq_ln_x}),
                 (\ref{eq:ln_eq_x}), (\ref{eq:lnq_e_x}) 
                 are equivalent to Eqs.\ (11a), (11b), (10a), (10b) 
                 of Ref.\ \protect\cite{BorgesCosta2021} respectively.
                 The notations introduced in Eqs.\ (\ref{eq:e_lnq_x}) 
                 and (\ref{eq:eq_ln_x}) 
                 differ from those used in \protect\cite{BorgesCosta2021}.}. 
Other generalizations exist in the literature, also referred to as $q$-numbers 
\cite{Haran2001,KacCheung2002}.
Generalized arithmetic operations follow from each of the $q$-numbers
and, consistently, there are various possibilities.
The present paper will only explore one possibility for $q$-numbers,
namely our Eq.\ (\ref{eq:e_lnq_x}), 
equivalent to the iel-number Eq.\ (11a) of \cite{BorgesCosta2021}.
For this choice, two algebras will be focused on here, namely,
\begin{equation}
\label{eq:operators-a}
\langle x \rangle_q \;\Circle^{q} \;\langle y \rangle_q 
  \equiv
  \langle\: x \: \circ \: y \,\rangle_q
\end{equation}
and
\begin{equation}
\label{eq:operators-b}
\langle\, x \: \Circle_{q} \:y \, \rangle_q 
  \equiv \,
  \langle x \rangle_q \circ \langle y \rangle_q \, ,
\end{equation}
the symbol $\circ$ representing any of the ordinary arithmetic operators
\mbox{$\circ \in \{+,-,\times,\slash\}$},
and $\Circle^q$ or $\Circle_q$ representing the corresponding 
generalized operators; 
naturally, 
$\Circle_{1} = \Circle^{1} =\circ$.
Possibility (\ref{eq:operators-a}) preserves prime number factorizability, 
while possibility (\ref{eq:operators-b}) does not
\footnote{
          The algebra corresponding to Eq.\ (\ref{eq:operators-a}) 
          is developed in Section \ref{factorizing-q-primes},
          and it was not addressed in Ref.\ \cite{BorgesCosta2021}.
          We use a superscript for its algebraic operators, $\Circle^{q}$,
          a notation that was not adopted in Ref.\ \cite{BorgesCosta2021}.
          The other algebra, corresponding to 
          Eq.\ (\protect\ref{eq:operators-b}),
          is presented in Section \protect\ref{non-factorizing-q-primes},
          and it corresponds to the oel-arithmetics addressed in Section III.D 
          of Ref.\ \protect\cite{BorgesCosta2021},
          where ${}_{\{q\}}\Circle$ is here noted $\Circle_q$.
         }.
%

\section{Preliminaries}
\label{preliminaries}

Let as remind, at this point, the so-called Basel problem, 
which focuses on the value of the series
\begin{equation}
 I_{2}=\sum_{n=1}^{\infty}\frac{1}{n^{2}} \, ,
\end{equation}
proposed in 1644 by Pietro Mengoli of Bologna, 
in contrast to the divergent harmonic series,
\begin{equation}
 I_{1}=\sum_{n=1}^{\infty}\frac{1}{n} \, .
\end{equation}
This problem was intensively studied by the Bernoulli brothers 
from Basel almost half a century later, 
and became known as the Basel Problem. 
They proved that $I_{1}$ is divergent, but $I_{2}$ is finite and smaller than 2,
but failed to obtain the exact value.

The Basel problem was solved by Euler in 1735. 
He also exhibited the connection with prime numbers, namely
$$
 I_{2}=\prod_{i=1}^\infty\frac{1}{1-p_{i}^{-2}}=\frac{1}{6}\pi^{2},
$$
where $p_{i}$ is the $i$-th prime numbers $(2,3,5,..)$, thus introducing 
for the first time the so-called Euler's product.
\begin{equation}
 I_{n}=\prod_{i=1}^\infty\frac{1}{1-p_{i}^{-n}}.
\end{equation}
In 1859, Riemann extended the domain of the exponent $n$ to complex numbers,
introducing the notation
\begin{align}
 \zeta\left(s\right)   & \equiv \sum_{n=1}^{\infty}\frac{1}{n^{s}}
\label{zeta-1}
\\
                       & =      \prod_{i=1}^\infty\frac{1}{1-p_{i}^{-s}},
 \label{zeta-2}
\end{align}
so that $\zeta\left(s\right)$ is often called Riemann's
zeta function (or Euler-Riemann zeta function).
By the way, Gauss made many important contributions to the field, 
especially the so-called prime number theorem, 
$\pi(N) \sim N/\ln N$ ($N \to \infty$),
where $\pi(N)$ is the number of primes up to the integer $N$.

For the oncoming discussion, it is useful to remind some properties 
within integers which are necessary to prove
\begin{equation}
 \sum_{n=1}^\infty\frac{1}{n^{s}} = 
 \prod_{i=1}^\infty\frac{1}{1-p_{i}^{-s}}. 
 \label{SUMProd}
\end{equation}
For this, let us consider a function $f(z)$, for $z = xy$,
which can be written as
\begin{equation}
f\left(  z\right)  =f\left(  x\right)  f\left(  y\right)  . \label{f=fxfy}%
\end{equation}
The function
$f\left(  z\right)  =\frac{1}{z^{s}}$
admits such a property. 
Indeed, the power law satisfies 
\begin{align}
z^{a+b}  &  =z^{a} z^{b},
\\
z^{-a}  &  =\frac{1}{z^{a}},\\
z^{a b}  &  =\left(  z^{a}\right)  ^{b}=\left(  z^{b}\right)  ^{a}.
\label{Powerlaws}%
\end{align}
Now, we know that any integer $n$ can be decomposed into an unique product 
of prime numbers%
\footnote{Let us focus on an interesting aspect of the primes. 
          Taking logarithm of
          Eq.\ (\ref{p-decompose}), for any positive integer 
          $n\in \mathbb{N} ,$ we get
          $\ln n  = \sum_{i}m\left(n\right)_{i}\ \ln p_{i}  
           = \sum_{i}m\left(n\right)_{i}\ \hat{e}_{i}$
          where $\hat{e}_{i} \equiv \ln p_{i}$.
          We can therefore consider the set 
          $\left\{\ln n,\ n\in \mathbb{N} \right\}$ 
          as a kind of infinite dimensional vector space, 
          whose basis vectors are 
          $\left\{\hat{e}_{i},i=1,..,\right\}.$ 
          However, rigorously, this is not a vector space, 
          since the coefficients of the linear combination, i.e., 
          the set of multiplicities, are only integers, 
          and not the set of real numbers.%
         },
\begin{equation}
 n=p_{_{1}}^{m_{1}} p_{2}^{m_{2}} ... p_{i}^{m_{i}} \cdots
 \label{p-decompose}
\end{equation}
where $m_{i}$ is the multiplicity of the prime $p_{i}$ in the product
($p_0=1,\, p_1=2,\, p_2=3,\, \dots$)
and they are \textit{uniquely} determined for a given 
$n\in \mathbb{N}.$ 
That is, the set
\begin{equation}
 \left\{  m_{1},m_{2},\cdots,m_{i},\cdots\right\}
\end{equation}
is determined for a given value of $n$ so we can use the notation
\begin{equation}
 \left\{m_{1}\left(n\right), 
        m_{2}\left(n\right), \cdots,
        m_{i}\left(n\right), \cdots
 \right\}\, .
\end{equation}
Note that to guarantee the uniqueness of the decomposition, the commutativity
$m \, n  = n \, m$
and the associativity
$l \, (m \, n) = (l \, m)  \, n$
of the product operation are essential.

Now, the sum over all positive integer values of a function 
$f(n)$ satisfying Eq.\ (\ref{f=fxfy}) can be written,
if it converges absolutely, as
\begin{equation}
 {\sum_{n=1}^\infty} f(n)
  = \sum_{n=1}^{\infty} \prod_{i}^\infty f(p_i^{m_{i}(n)}) \, ,
 \label{Euler_factorial}
\end{equation}

This comes from the fact that the direct product of the sets
$\left\{p_{i},p_{i}^{2},p_{i}^{3},\cdots,p_{i}^{m},\cdots\right\}$ 
for all primes
is equal to the set of positive natural numbers 
$\mathbb{N}^{+} \equiv \{1, 2, 3, \cdots \}$,
\begin{equation}
 \prod_{i=0}^\infty 
\scalebox{1.5}{$\otimes$} 
 \left\{p_{i},p_{i}^{2},p_{i}^{3}, \cdots,p_{i}^{m} ,\cdots\right\} 
 = \mathbb{N}^+,
 \label{SetPrimes}
\end{equation}
where \scalebox{1.5}{$\otimes$} denotes the direct product. 
The meaning of Eq.\ (\ref{Euler_factorial}) can be seen more intuitively here 
below as
\begin{align}
&
\begin{array}[c]{ccccccc}
\;\;    \{1 & 2 & 2^{2} & 2^{3} & \cdots & 2^{m}\cdots & \}
 \\
  \text{\scalebox{1.5}{$\otimes$}}\{1 & 3 & 3^{2} & 3^{3} & \cdots & 3^{m}\cdots & \}
 \\
\text{\scalebox{1.5}{$\otimes$}}\{1 & \vdots & \vdots & \vdots & \cdots & \vdots & \}
 \\
\text{\scalebox{1.5}{$\otimes$}}\{1 & p_{i} & p_{i}^{2} & p_{i}^{3} & \cdots & p_{i}^{m} & \}
 \\
 \vdots & \vdots & \vdots & \vdots & \vdots & \vdots & \vdots
\end{array}
\nonumber \\
&
= \left\{  1,1\cdot 2,\cdots,p_{_{1}}^{m_{1}}\cdot p_{2}^{m_{2}}%
\cdot ...\cdot p_{i}^{m_{i}}\cdot ,\dots\right\} 
= \mathbb{N} \, .
\label{prime-set}
\end{align}

Note that the multiplicity $m_{i}(n)  $ of the $i$-th prime
$p_{i}$ takes all integer values if $n$ runs over all natural numbers.
Consistently, for a given $p_{i},$ the value of $m_{i}(n)$ runs
over all nonnegative integers. That is, $m_{i}(n)$ and $i$ can run over
all nonnegative integers independently. 
Therefore, by using $l\,(m+n) = l\,m + l\,n$,
we can exchange the order of the sum and the product
in Eq.\ (\ref{Euler_factorial}):
\begin{equation}
\sum_{n=1}^\infty f\left(  n\right)  
= \prod_{i=1}^\infty\sum_{m=0}^{\infty}f\left(p_{i}^{m}\right)  
\label{Euler Sum}
\end{equation}
if $f(x)$ satisfies the property indicated in Eq.\ (\ref{f=fxfy}). 

Now, let us take the case of the Euler product,%
\[
f\left(  n\right)  =\frac{1}{n^{s}}.
\]
Then the summation in $m$ of Eq.(\ref{Euler Sum}) for $p_{i}$ term is
\begin{equation}
S\equiv\sum_{m=0}^\infty \left(  \left(  \frac{1}{p_{i}}\right)  ^{s}\right)
^{m}=1+\frac{1}{p_{i}^{s}}+\left(  \frac{1}{p_{i}^{s}}\right)  ^{2}
+\cdots+\left(  \frac{1}{p_{i}^{s}}\right)  ^{k}+\cdots
\end{equation}
which is a geometric series. 
By writing $r \equiv 1/p_{i}^{s}$,
we obtain
\begin{equation}
S \equiv 1 + r + r^{2} + \cdots = \frac{1}{1-r} 
= \frac{1}{1-\left(  1/p_{i}\right)  ^{s}} \, .
\end{equation}
Finally, we have the Euler product form as%
\begin{align}
 \sum_{n=1}^\infty \frac{1}{n^{s}}  
 = \prod_{i=1}^{\infty}\sum_{m=0}^{\infty}\left( \frac{1}{p_{i}^{s}}\right)^{m}
 = \prod_{i=1}^{\infty}\left(\frac{1}{1-\left(  1/p_{i}\right)  ^{s}}\right) 
 \, .
\end{align}
This connection is known as Euler's product form. 
Riemann extended the domain of $s$ of this function into the complex
plane $z$, being since then frequently referred to as 
the Riemann zeta function $\zeta(z)$.

\section{$q$-integers}

We focus on the $q$-generalized numbers $\langle n \rangle_q$,
$n \in \mathbb{N}$, given by Eq.\ (\ref{eq:e_lnq_x}),
to see how far we can conserve the essential concept of prime numbers 
in the set of $q$-numbers.
This particular $q$-number $\langle n \rangle_q$ satisfies
\begin{equation}
 \ln\, \langle n \rangle_q = \ln_q n.
 \label{eq:qlog_x.eq.log_xq}
\end{equation}

Let us consider the set of nonnegative integers $\mathbb{N}$.
Eq.\  (\ref{eq:e_lnq_x}) defines a mapping from 
$\mathbb{N}$ to $\mathbb{N}_{q}$, where $\mathbb{N}_{q}$ denotes
the set of all positive $q$-integers. 
Note that 
$\lim_{x\to 1} \langle x \rangle_q \to 1$
but, for other general natural numbers $n$, its $q$-partner
\begin{equation}
 \langle n \rangle_q =\exp\left\{  \frac{n^{1-q}-1}{1-q}\right\}
 \label{q-integer}
\end{equation}
is not an integer, in general. 
This point is crucial for the question of factorization.
The inverse mapping is then (see Eq.\ (\ref{eq:eq_ln_x}))
\begin{equation}
  n = \exp_{q}\left\{  \ln\left( \langle n \rangle_q \right ) \right\} \\
    = {}_q\langle\, \langle n \rangle_q \, \rangle.
 \label{q_n}%
\end{equation}
The $q$-mapping and its inverse
$\mathbb{N} \leftrightarrows \mathbb{N}_{q}$
are one-to-one
[$ \exp_{q}(\ln_{q} x) = 
 \ln_{q}(\exp_{q}(x)) = x, \,\forall x\in \mathbb{R}_+$].
However, it is clear that, for $q\neq1$,
$\mathbb{N} \neq \mathbb{N}_{q}$.

\begin{figure}[h!]
\centering
\includegraphics[width=9.5cm]{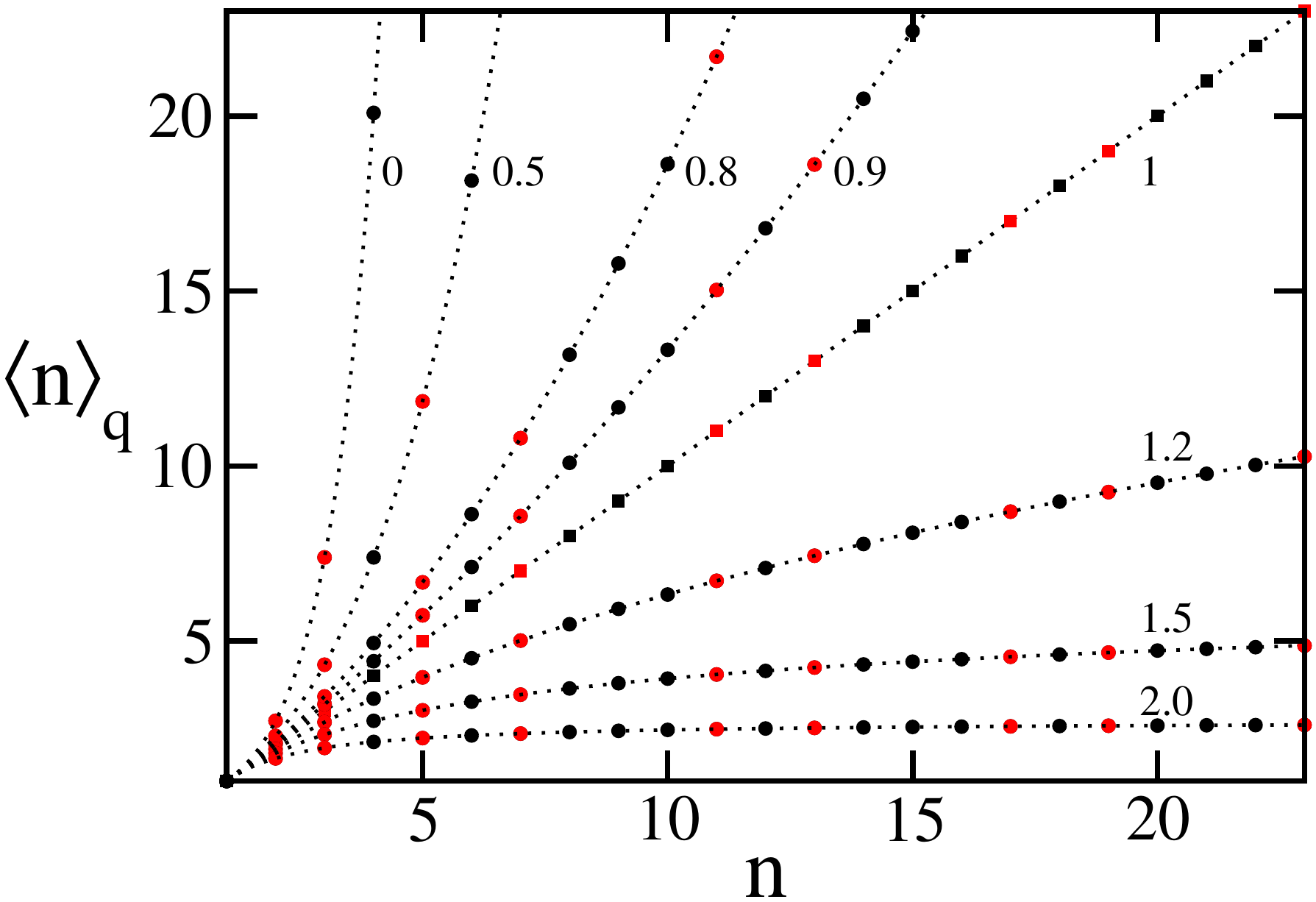}
\caption{\label{fig1} 
         \small
         $q$-natural numbers for typical values of $q$. 
         The $q$-prime numbers are indicated in red. 
         The dotted lines correspond to real values for the abscissa.
        }
\end{figure}
For $q \in (1,\infty)$, the $q$-natural numbers satisfy 
$ \lim_{n\to\infty} \langle n\rangle_q =e^{1/(q-1)} >1\,$.
For $q \in (-\infty,1]$, the $q$-natural numbers satisfy
$ \lim_{n\to\infty} \langle n\rangle_q =\infty  \,$.
There are infinitely many 
$q$-prime numbers for $q \in (-\infty,\infty)$.  
See Fig.\ \ref{fig1}.

\section{\label{factorizing-q-primes}
         Algebra preserving factorizability of $q$-integer numbers
         in $q$-prime numbers}

In order to introduce the concept of $q$-primes, we should keep the
factorization concept in $\mathbb{N}_{q}^{+}$. 
Thus we must first define the product operation $\mathbb{N}_{q}$, 
and should be kept invariant under the $q$-mapping from $\mathbb{N}^{+}$. 
Following Eq.\ (\ref{eq:operators-a}), let 
$\langle m \rangle_q \otimes^{q} \langle n \rangle_q$
be such a product%
\footnote{Note that here $\otimes^{q}$ is not the direct product of sets, 
          but this specific $q$-generalized product of $q$-numbers.}           
between any pair of $q$-integers 
$ \langle m \rangle_q, \, \langle n \rangle_q \in \mathbb{N}_{q}$,
whose result give the $q$-transform of $m\,n\ \in \mathbb{N} $, 
i.e., obeys the following factorizability:
\begin{equation}
 \langle m \rangle_q \otimes^q \langle n \rangle_q = \langle m n \rangle_q.
 \label{eq:q-product-a}
\end{equation}
We also define the following generalized summation operator 
$\oplus^q$ in $\mathbb{N}_{q}$ as
\begin{equation}
 \langle m \rangle_q \oplus^{q} \langle n \rangle_q = \langle m + n \rangle_q.
 \label{eq:q-sum-a}
\end{equation}
The above definitions correspond to those introduced 
in Ref.\ \cite{Lobao-etal2009}
[its Eqs.\ (5) and (6)]%
\footnote{The operations defined in \protect\cite{Lobao-etal2009} 
          follows the same structure of Eq.\ (\protect\ref{eq:q-product-a}), 
           or, more generally, Eq.\ (\protect\ref{eq:operators-a}),
          but with a different deformed number, namely the Heine number.
         }.

In order that the above operations are meaningful in $\mathbb{N}_{q}^{+}$, 
it is necessary that they are closed operations in $\mathbb{N}_{q}$. 

The following definitions satisfy
properties (\ref{eq:q-product-a}) and (\ref{eq:q-sum-a}): 
\begin{equation}
 \langle m \rangle_q  \otimes^q  \langle n \rangle_q
 = \langle m \rangle_q \; \langle n \rangle_q \;
    e^{(1-q) \left(\ln \langle m \rangle_q \right) 
             \left( \ln \langle n \rangle_q\right)},
 \label{eq:q-product-a-def}
\end{equation}
\begin{equation}
 \langle m \rangle_q  \oplus^q  \langle n \rangle_q
  = \exp\left\{ \frac{(m+n)^{1-q}-1}
                {1-q} 
        \right\}.
 \label{eq:q-sum-a-def}
\end{equation}
\\
Eq.\ (\ref{eq:q-product-a-def}) can be rearranged as
\begin{equation}
 \langle m \rangle_q  \otimes^q  \langle n \rangle_q
 = \langle m \rangle_q \; \langle n \rangle_q \;
   e^{\ln \langle m \rangle_q \,  \oplus^q \, \ln \langle n \rangle_q
       \, - \, \ln \langle m \rangle_q \, - \, \ln \langle n \rangle_q }.
\end{equation}
The result of these definitions is a $q$-integer by construction.

By construction, it is clear that this definition of the $q$-product
$\otimes^{q}$ as an operator in $\mathbb{N}_{q}$, 
conserves the factorization property under the $q$-transform 
$
 n  = k\,m \in \mathbb{N} \Longleftrightarrow\ 
 \langle n \rangle_q = \langle k \rangle_q \otimes^q \langle m \rangle_q 
 \in \mathbb{N}_{q}. 
$
In addition, also by construction, the operators, $\otimes^{q}$ 
and $\oplus^{q}$ satisfy the following basic properties of algebras, 
valid in $\mathbb{N}_{q}$: 

\begin{itemize}
\item Closedness of the operation:
 \begin{equation}
  \forall \langle m \rangle_q, \langle n \rangle_q \in \mathbb{N}_{q}^{+}
  \Longrightarrow\
  \langle m \rangle_q \oplus^q \langle n \rangle_q
  \in \mathbb{N}_{q}^{+}, 
  \, 
  \langle m \rangle_q \otimes^q \langle n \rangle_q
  \in \mathbb{N}_{q}^{+}, 
  \label{closure-q_sup}
 \end{equation}

\item Commutativity
\begin{align}
 \langle m \rangle_q \oplus^q \langle n \rangle_q
 &=
 \langle n \rangle_q \oplus^q \langle m \rangle_q,
 \label{commutativity-plus-sup_q}
\\
 \langle m \rangle_q \otimes^q \langle n \rangle_q
 &=
 \langle n \rangle_q \otimes^q \langle m \rangle_q,
\label{commutativity-q_sup-product}
\end{align}

\item Associativity
 \begin{align}
  \langle k \rangle_q \oplus^q
  \bigl(\langle m \rangle_q \oplus^q \langle n \rangle_q\bigr)
  &=
  \bigl(\langle k \rangle_q \oplus^q \langle m \rangle_q\bigr)
  \oplus^q \langle n \rangle_q,
 \label{associativity-q_sup-sum}
\\
  \langle k \rangle_q \otimes^q
  \bigl(\langle m \rangle_q \otimes^q \langle n \rangle_q\bigr)
  &=
  \bigl(\langle k \rangle_q \otimes^q \langle m \rangle_q\bigr)
  \otimes^q \langle n \rangle_q,
 \label{associativity-q_sup-prod}
\end{align}

\item Distributivity of the product $\otimes^q$ 
      with regard to the sum $\oplus^q$
 \begin{equation}
  \langle k \rangle_q \; \otimes^q \,
  \bigl( \langle m \rangle_q \, \oplus^q  \, \langle n \rangle_q \bigr)
  =
  \bigl( \langle k \rangle_q \; \otimes^q \, \langle m \rangle_q \bigr)
  \;
  \oplus^q
  \;
  \bigl( \langle k \rangle_q \; \otimes^q \, \langle n \rangle_q \bigr).
 \label{distributivity-q_sup}
 \end{equation}

\item Neutral element of the $q$-addition
\begin{equation}
 \langle m \rangle_q \oplus^q \langle 0 \rangle_q 
 = \langle m \rangle_q \oplus^q 0
 = \langle m \rangle_q\;\;\;(q \ge 1).
 \label{neutral-q_sup-sum}
\end{equation}

\item Neutral element of the $q$-product
\begin{equation}
 \langle m \rangle_q \otimes^q \langle 1 \rangle_q 
 = \langle m \rangle_q \otimes^q 1 
 = \langle m \rangle_q \,.
 \label{neutral-q_sup-prod}
\end{equation}

\end{itemize}

We show that these properties are essential to keep the nature of prime numbers 
for $q$-transformed integers, corresponding to Eq.\  (\ref{p-decompose}) 
in $\mathbb{N}$.

We also define the following operations:
\begin{equation}
  \langle x \rangle_q \oslash^q \langle y \rangle_q = \langle x/y \rangle_q
\end{equation}
and
\begin{equation}
  \langle x \rangle_q \ominus^q \langle y \rangle_q = \langle x-y \rangle_q
\end{equation}

Besides, for any $a \in \mathbb{R}$, 
we can define a $a$-power of a $q$-number by
\begin{equation}
 \langle x \rangle_q \owedge^q a \equiv \Bigl\langle x^a \Bigr\rangle_q
 \label{PowerInteger}
\end{equation}
from what follows
\begin{align}
 \langle x \rangle_q \owedge^q (a+b) 
  &=
     \Bigl\langle x^{a+b} \Bigr\rangle_q
 \nonumber\\
  &= \langle x^a \rangle_q \, \otimes^q \, \langle x^b \rangle_q
 \nonumber\\
  &= \Bigl(\langle x \rangle_q \owedge^q a\Bigr) 
     \, \otimes^q \, 
     \Bigl(\langle x \rangle_q \owedge^q b\Bigr),
 \label{power(a+b)}
\end{align}
\begin{align}
 \langle x \rangle_q \owedge^q (ab) 
 &= \Bigl\langle x^{ab} \Bigr\rangle_q
\nonumber\\
 &= \Bigl\langle x^a \Bigr\rangle_q \owedge^q b
\nonumber\\
 &= \Bigl( \; \langle x \rangle_q \owedge^q a \; \Bigr) \owedge^q b
\nonumber\\
 &= \Bigl( \; \langle x \rangle_q \owedge^q b \; \Bigr) \owedge^q a,
 \label{power(axb)}
\end{align}
and
\begin{align}
  \Bigl(\langle x \rangle_q \owedge^q s\Bigr) 
     \, \otimes^q \, 
  \Bigl(\langle x \rangle_q \owedge^q t\Bigr) 
 &= \Bigl(\langle x \rangle_q\Bigr)^{\langle s \rangle_q 
                                      \, \oplus^q \, 
                                      \langle t \rangle_q
                                     }.
\end{align}
As we see, for fixed $q$, this algebra is isomorphic to the standard algebra.

\noindent
{\bf Definition 1\ }
{\it 
A $q$-integer 
$\langle n \rangle_q$
(i.e., an element of $ \mathbb{N}_{q}^{+}$) 
is called ``$q$-prime'' and written as 
$\langle p \rangle_q$
if it can not be written as a $q$-factorized form in terms of two smaller 
$q$-integers as
\begin{equation}
 \langle n \rangle_q = \langle m \rangle_q \otimes^q \langle k \rangle_q,
\end{equation}
with
$\langle m \rangle_q, \langle k \rangle_q \in \mathbb{N}_q^+$,
except for the trivial factorization case, i.e., either one of 
$\langle m \rangle_q$ or $\langle k \rangle_q$
is the unity, $\langle 1 \rangle_q$.
}

It is evident from the definition of $q$-integers 
that all the set of $q$-primes 
$\left\{ \langle p_i \rangle_{q}, i=1,\dots \right\}$ are $q$-partners
of the prime numbers in $ \mathbb{N}^{+}$, $\left\{  p_{i},i=1,...\right\}$.
For any integer $n\in \mathbb{N}^{+}$ which is not a prime, 
there exists the non-trivial factors, $k, m \in \mathbb{N}^{+}$ such that
$ n=k \times m$.
But from the definition of the $q$-product, Eq.\ (\ref{eq:q-product-a}),
$\langle n \rangle_q
 = \langle k m \rangle_q
 = \langle k \rangle_q \, \otimes^q \, \langle m \rangle_q$,
showing that $\langle n \rangle_q$ has a non-trivial $q$-factorization 
and it is not a $q$-prime. 
As we mentioned, the $q$-correspondence between the two sets, 
$ \mathbb{N}^{+}$ and $ \mathbb{N}_{q}^{+}$ is one-to-one, 
$q$-primes in $ \mathbb{N}_{q}^{+}$ are $q$-transformations 
of the normal primes in $ \mathbb{N}^{+}$.

The basic property of primes is that any natural number
$n \ge 2$ can be written uniquely as the products of primes 
$\left\{p_{i}\right\}$ 
as
\begin{equation}
 n=\prod_{i}p_{i}^{m_{i}(n)}, 
 \label{decomposition}
\end{equation}
where
$$ \left\{p_{1},p_{2},p_{3},p_{4},p_{5}\ \cdots\ \ \right\} 
   = \left\{2,3,5,7,\ \cdots\right\}$$ 
is the set of primes and for a given $n$, 
$$ \left\{m_{1}(n), m_{2}(n), \cdots\right\} $$
is the set of multiplicities $m_{i}(n)$ of the prime $p_{i}$ 
is uniquely determined for given $n$.

Now, for $q$-integers, 
$\langle n \rangle_q \in \mathbb{N}_{q}^{+}$, 
the corresponding decomposition property is valid in $\mathbb{N}_{q}^{+}$ 
in terms $q$-integers with $q$-products, satisfying the properties
of commutativity and distributivity. 
By construction, we can write for any $\langle n \rangle_q$
the $q$-prime decomposition as
\begin{equation}
 \langle n \rangle_q = {\displaystyle\prod_i}^q
                         \;\; 
                         \Bigl(
                               \langle p_i \rangle_q \owedge^q m_i(n)
                         \Bigr)
\end{equation}
for the sake of isomorphism of the product operations in 
$\mathbb{N}$ and $\mathbb{N}_{q}$.
The symbol ${\prod}^q$ represents the generalized product (\ref{eq:q-product-a})
of a number of terms,
$\overset{\!\!n}{{\underset{\!\!i}{\sideset{}{^q}\prod}}}
  x_i \equiv x_1 \otimes^q x_2 \otimes^q \cdots \otimes^q x_n$.

Since we have the isomorphisms of the operations of sum and product in 
$\mathbb{N}^{+}$ and in $\mathbb{N}_{q}^{+}$, 
we can write down the $q$-version of the Euler product:
\begin{equation}
 {\sum_{\langle n \rangle_q \in \mathbb{N}_{q}^{+}}}^{\!\!\!\!\!q}
   \; \left(  \frac{1}{\langle n \rangle_q} \right)^{m \, \oplus^{q}\, s} 
  = 
 {\sum_{\langle n \rangle_q \in \mathbb{N}_{q}^{+}}}^{\!\!\!\!\!q}
 \;\prod_{i}
 \;\Bigl(\langle p_i \rangle_q\Bigr)^{m_i(n)},
\end{equation}
with the symbol ${\sum}^q$ representing the generalized summation
according to Eq.\ (\ref{eq:q-sum-a}).

Now, the multiplicity $m_{i}(n)$ of the $i$-th prime $p_{i}$
should take all integer values if $n$ runs over all the integers. 
Inversely, for a given $p_{i}$, the value of $m_{i}(n)$ runs over all integers. 
That is, $m_{i}(n)$ and $i$ can run over all integers independently. 
This comes from the fact that the direct product of
the sets, $\left\{p_{i},p_{i}^{2},p_{i}^{3},\cdots,p_{i}^{m},\cdots\right\}$, 
for every prime is equal to the set of natural numbers $\mathbb{N}^{+}$
(see Eq.\ (\ref{SetPrimes})).
We can therefore see that exchange the order of the sum and the product%
\footnote{
          Due to the distributivity of the ordinary multiplication 
          with regard to the ordinary addition,
          $l\,(m+n) = l\,m + l\,n$.
         } 
in Eq.\ (\ref{Euler_factorial}) as
\begin{equation}
 \sum_{n\in \mathbb{N}^{+}}f(n)  
 = \sum_{n=1}^{\infty} \prod_{i} f\left(p_{i}^{m_{i}(n)}\right)  
 = \prod_{i=1} \sum_{m=0}^{\infty} f\left(p_{i}^{m}\right)
\end{equation}
if $f\left(  x\right)  $ satisfies the property Eq.\ (\ref{f=fxfy}).

\begin{figure}[h!] 
\centering
 \includegraphics[width=7.9cm]{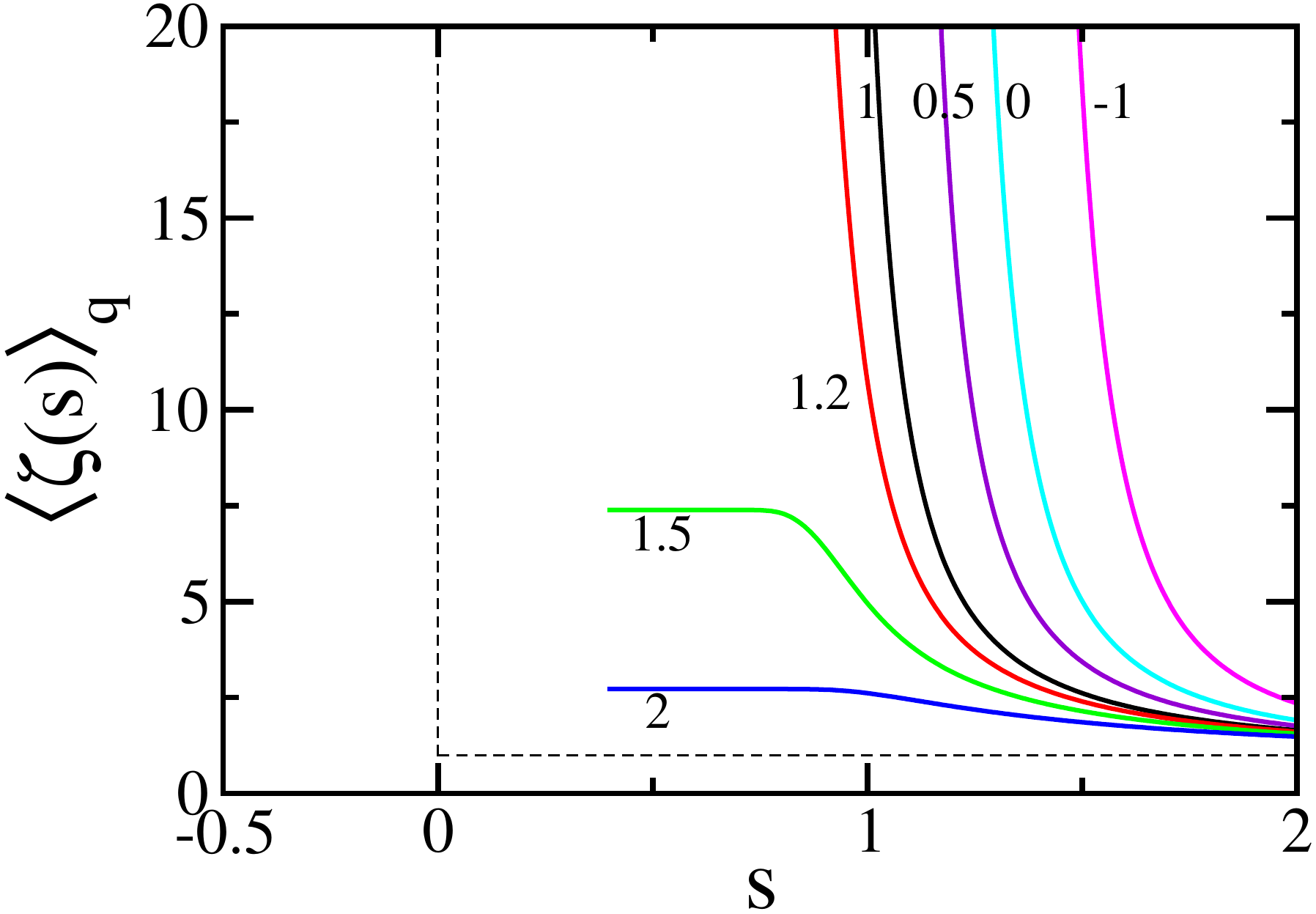}
\caption{\label{fig:q-zeta1} 
         \small
         $\langle \zeta (s)\rangle_q=\zeta_q(s)$ 
         for typical values of $q$. 
         These values have been calculated from 
         Eq.\ (\protect\ref{eq:zeta_q.eq.q-zeta_1})
         with the first $10^5$ prime numbers. 
         The $q$-number $\langle n \rangle_q$, Eq.\ (\protect\ref{q-integer}),
         with $q>1$ presents an upper limit, 
         $\lim_{n \to \infty} \langle n \rangle_{q>1} = e^{1/(q-1)}$,
         and this prevents the divergence of $\zeta_q(s)$ for $s<1$.
        }
\end{figure}

The definitions of the $q$-algebra introduced here is sufficient to write down
the $q$-version of the Euler sum and the Euler product in $q$-representation
formally as
\begin{align}
 \zeta_q (s)   
  &= \underset{\forall \langle n \rangle_q \in \mathbb{N}_{q}}
              {\sum}^{\!\!\!\!\!\!q} \;\;
                       \Bigl\langle  \;
                                   \langle 1 \rangle_q 
                                   \, \oslash^q \, 
                                   \left(\langle n \rangle_q \owedge^q s \right)
                                   \;
                       \Bigr\rangle_q
\label{eq:zeta-q-sums}
\\
 &= {\underset{i \in \mathbb{N}}{\sideset{}{^q}\prod}} \;
    \Bigl\langle \;
                 \langle 1 \rangle_q
                 \, \oslash^q \, 
                 \langle \, 1 - 1/p_i^s \, \rangle_q
    \Bigr\rangle_q,
 \label{eq:zeta-q-a}
\end{align}
and 
\begin{equation}
 \label{eq:zeta_q.eq.q-zeta_1}
 \zeta_q(s) \equiv \langle \zeta(s) \rangle_q
\end{equation}
(with $\zeta_1(s) \equiv \zeta(s)$).
In other words, the Euler product form is preserved for all values of $q$.
Fig.\ \ref{fig:q-zeta1} depicts $\zeta_q(s)$ for different values of $q$.

\begin{figure}
    \centering
    \includegraphics[width=0.30\textwidth]
     {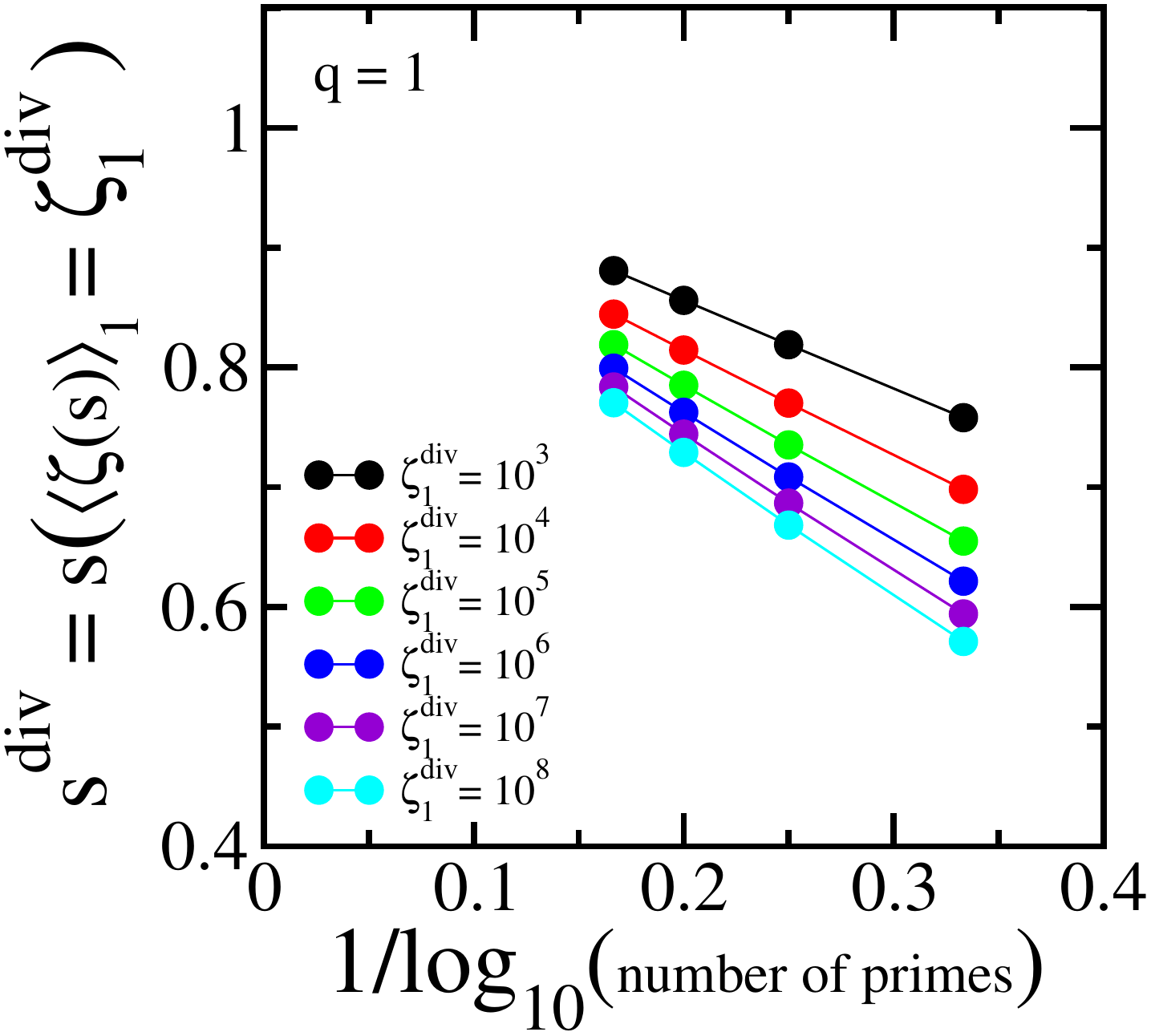}
    \includegraphics[width=0.30\textwidth]
     {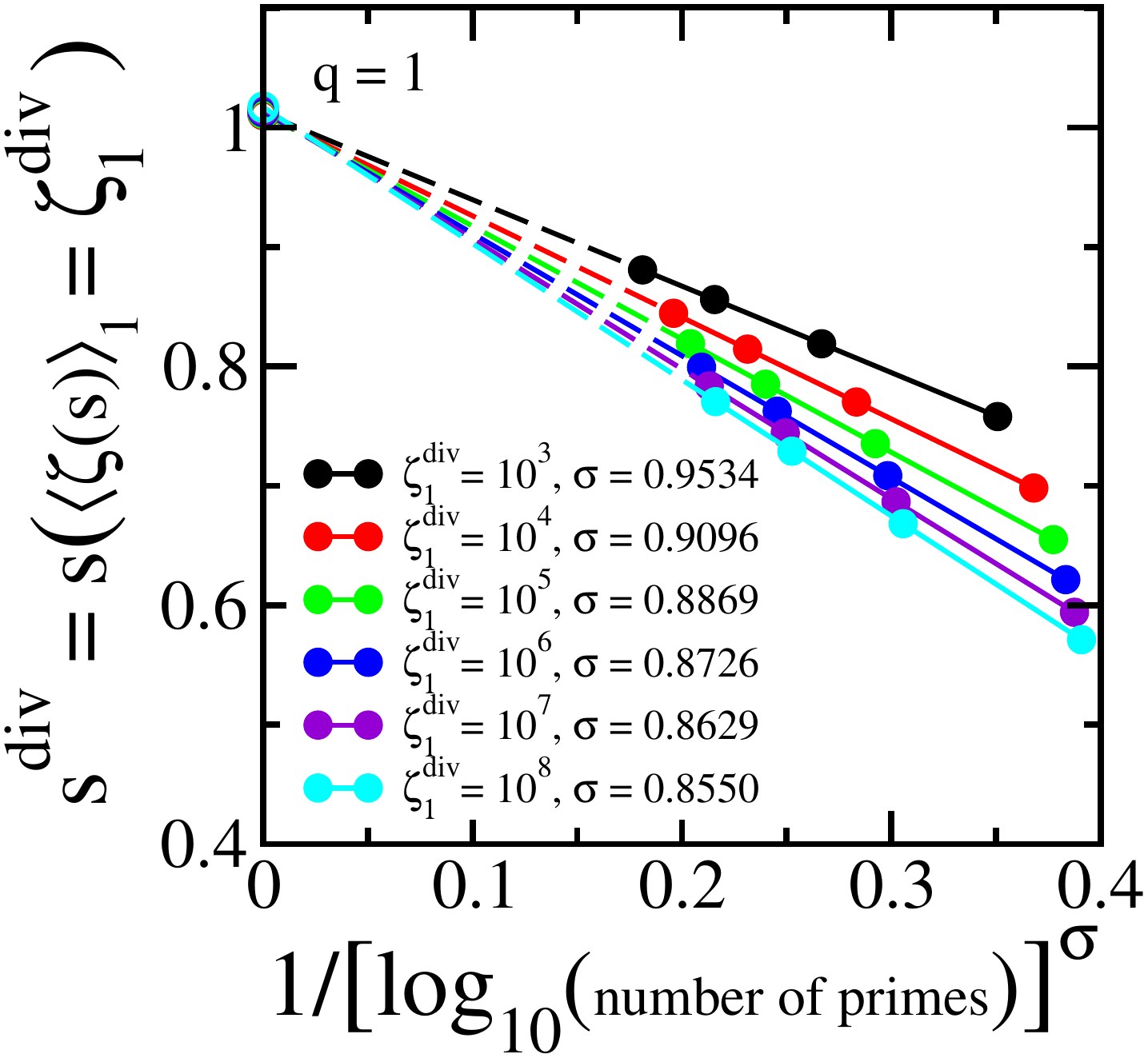}
    \includegraphics[width=0.30\textwidth]
     {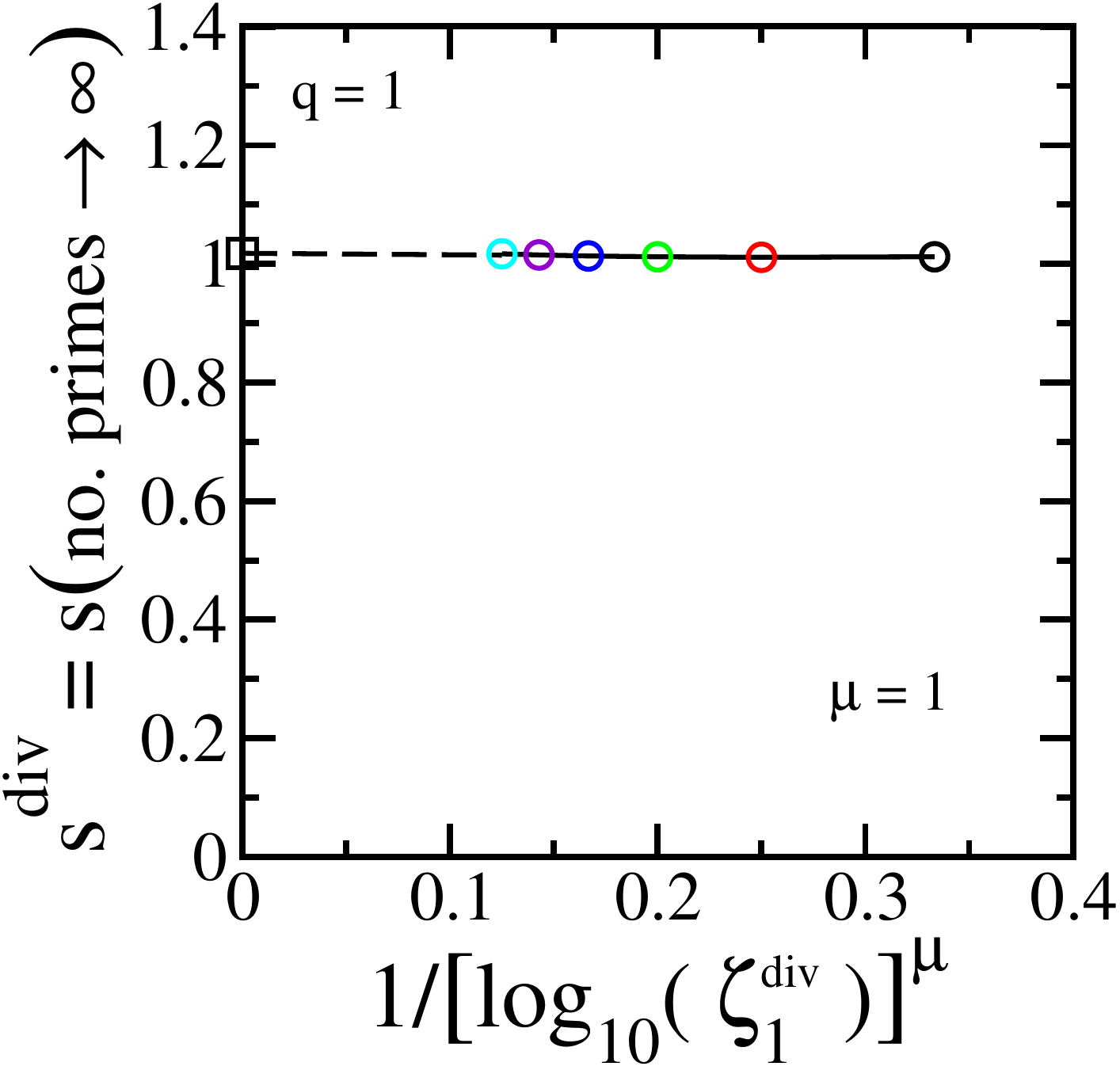}

    \includegraphics[width=0.30\textwidth]
     {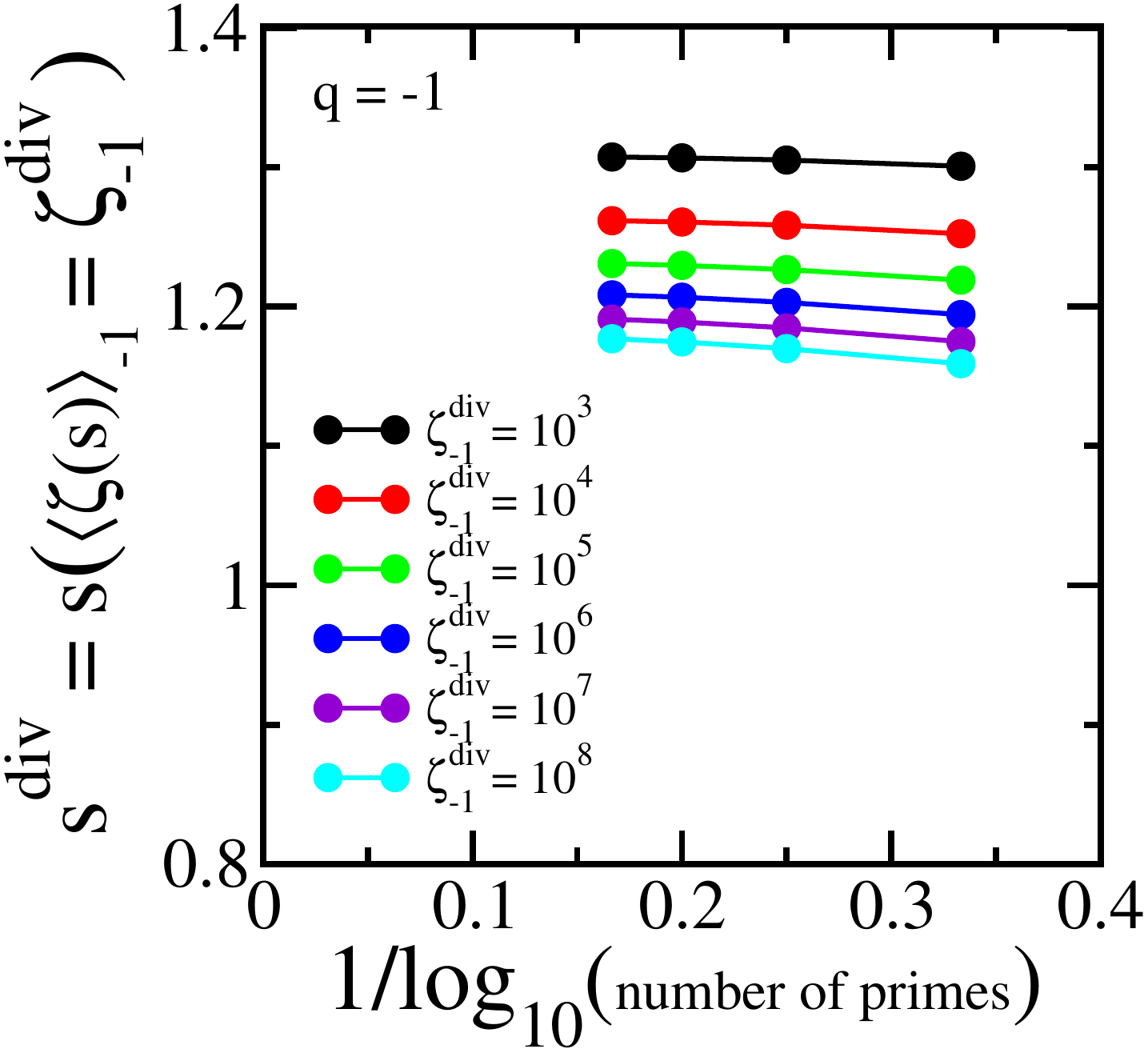}
    \includegraphics[width=0.30\textwidth]
     {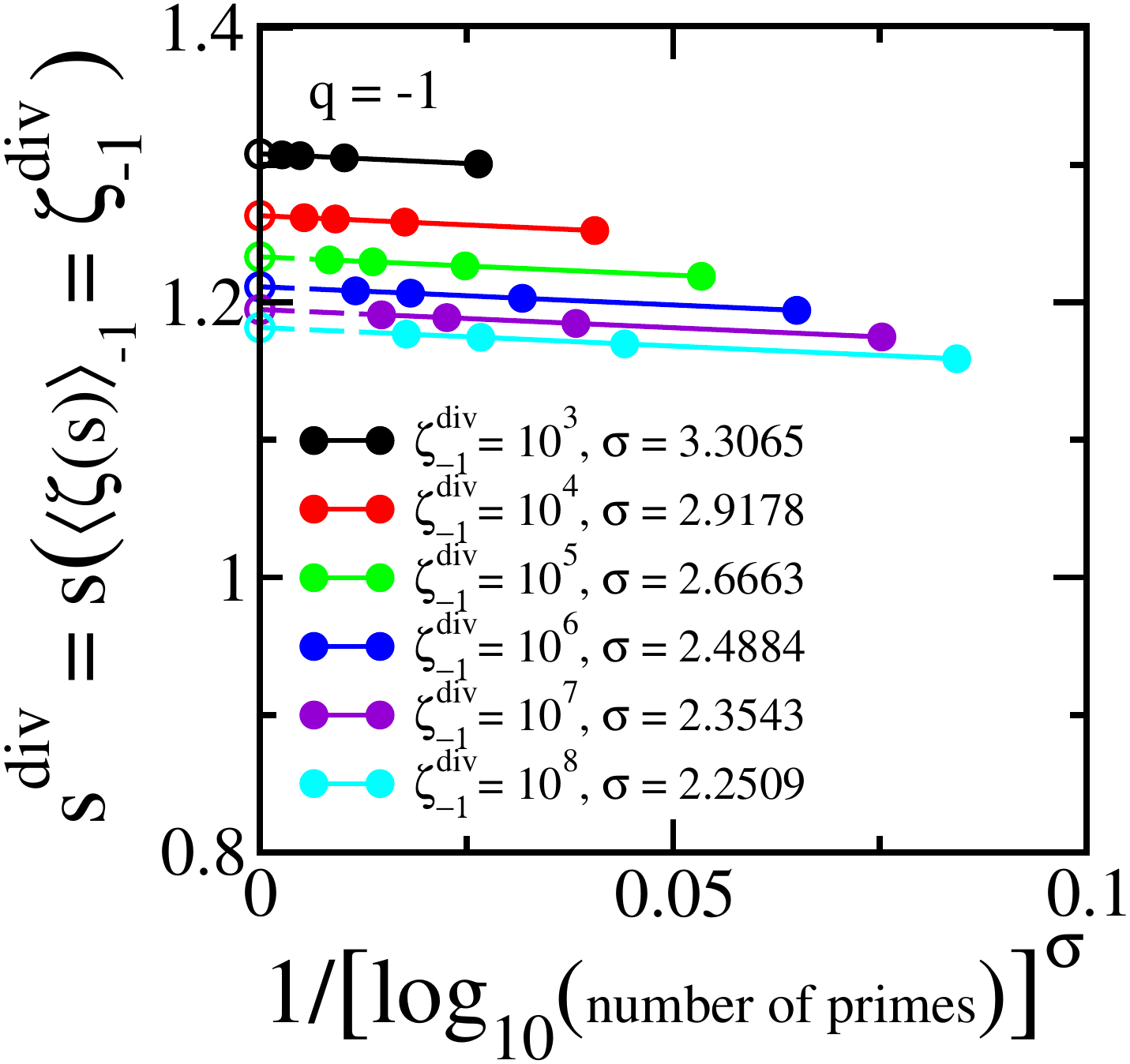}
    \includegraphics[width=0.30\textwidth]
     {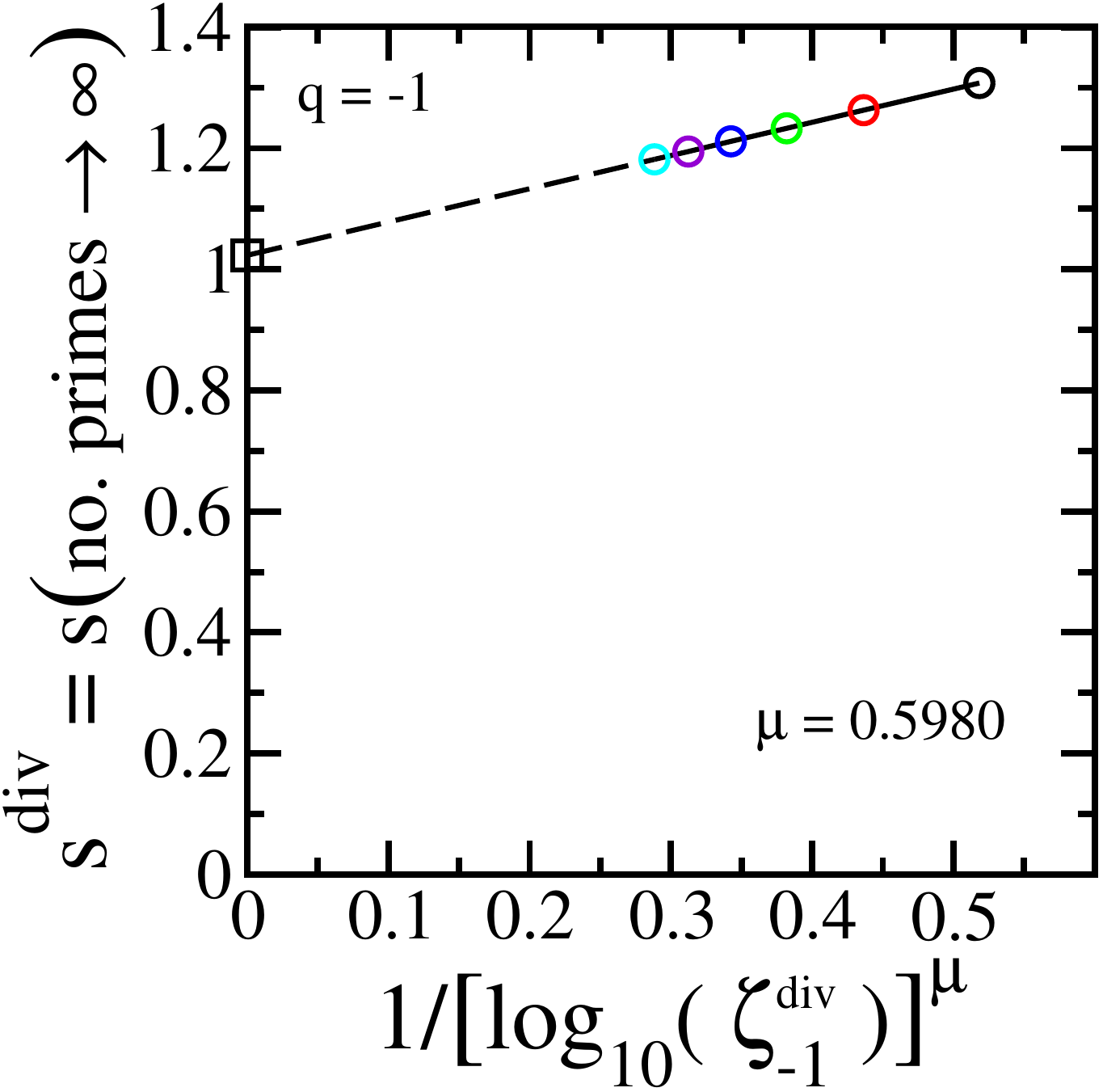}
\caption{\label{fig:zeta1} 
         \small
         {\it Top left panel}:
         Dependence of  
         $s^{\text{div}}=s(\zeta_1^{\text{div}})$ with 
         $1/\log_{10}(\text{number of primes})$
         in Eq.\ (\protect\ref{zeta-2});
         $\zeta_1^{\text{div}}$ is a proxy for the divergence of $\zeta_1(s)$. 
         {\it Top middle panel}:
         Power-law rescaling of 
         $s^{\text{div}}=s(\zeta_1^{\text{div}})$ with 
         $1/\left[\log_{10}(\text{number of primes})\right]^\sigma$.
         The curves point towards $s^{\text{div}}$ 
         with infinite number of primes (open circles);
         {\it Top right panel}:
         Extrapolated values of Fig.\ \ref{fig:zeta1}(top middle)
         linearly rescaled with 
         $1/\left[\log_{10}(\zeta_q^{\text{div}})\right]^\mu$
         point towards the analytically exact value
         $s^{\text{div}(\infty)}=1$ (open square)
         ($ \lim_{\zeta_1^{\text{div}}\to\infty} 
            \lim_{\text{number of primes}\to\infty} s^{\text{div}} = 1$)
         within a numerical error less than 2\% for $q=1$.
         The colors of the open circles refer to the values of
         $\zeta_1^{\text{div}}$ identified in Fig.\ \ref{fig:zeta1}(top middle).
         {\it Bottom panels}: the same as top panels for $q=-1$,
         Eq.\ (\protect\ref{eq:zeta_q.eq.q-zeta_1}). 
        }
\end{figure}

We numerically identify the location of the divergence
by fixing an arbitrary value of $\zeta_q(s)$ noted as $\zeta_q^{\text{div}}$ 
and identify the corresponding the value of $s$ (noted as $s^{\text{div}}$)
with increasing number of primes. 
The procedure is repeated with increasing values of $\zeta_q^{\text{div}}$.
The three top panels of Fig.\ \ref{fig:zeta1} illustrate the procedure 
for $q=1$, and the three bottom panels for $q=-1$.
Each curve in Fig.\ \ref{fig:zeta1}(top left)
displays the  value of $s^{\text{div}}$ for which 
$\zeta_q(s^{\text{div}}) \equiv \zeta_q^{\text{div}} 
 = \langle \zeta_1(s^{\text{div}}) \rangle_q
 = 10^3, 10^4, \dots, 10^8$
with increasing numbers of primes 
($10^3, 10^4, \dots, 10^6$ primes, shown with solid circles).
The representation with $1/\log_{10}(\text{number of primes})$ in the abscissa
is not a straight line, and we empirically found that introducing a power
$\sigma$ ($\sigma$ depends on $\zeta_q^{\text{div}}$)
as shown in Fig.\ \ref{fig:zeta1}(top middle),
straight lines emerge, which can be extrapolated 
(dashed lines)
to infinite number of primes --- the open circles at the ordinate axis.
These extrapolations correspond to infinite number of primes,
but, nevertheless, the values of $\zeta_q^{\text{div}}$ are still finite.
Finally, the limit $\zeta_q^{\text{div}} \to \infty$
is achieved as illustrated in Fig.\ \ref{fig:zeta1}(top right).
The open circles at the ordinate axis of Fig.\ \ref{fig:zeta1}(top middle) 
are represented in Fig.\ \ref{fig:zeta1}(top right) for each value of 
$\zeta_q^{\text{div}}$, identified with their respective colors,
with the change of variables
$1/[\log_{10}(\zeta_q^{\text{div}})]^\mu$
in the abscissa.
$\mu$ is an empirical power that transforms the curves into straight lines
($\mu$ depends on $q$; $\mu(q=1)=1$). A final extrapolation is then allowed,
identifying the location of the divergence (open square).
The difference 
$ \lim_{\zeta_1^{\text{div}}\to\infty} 
  \lim_{\text{number of primes}\to\infty} s^{\text{div}} - 1$
characterizes the numerical error. 
We use the same procedure for $q < 1$, and the three bottom panels 
illustrate the case $q=-1$.

\begin{figure}
    \centering
    \includegraphics[width=0.30\textwidth]
     {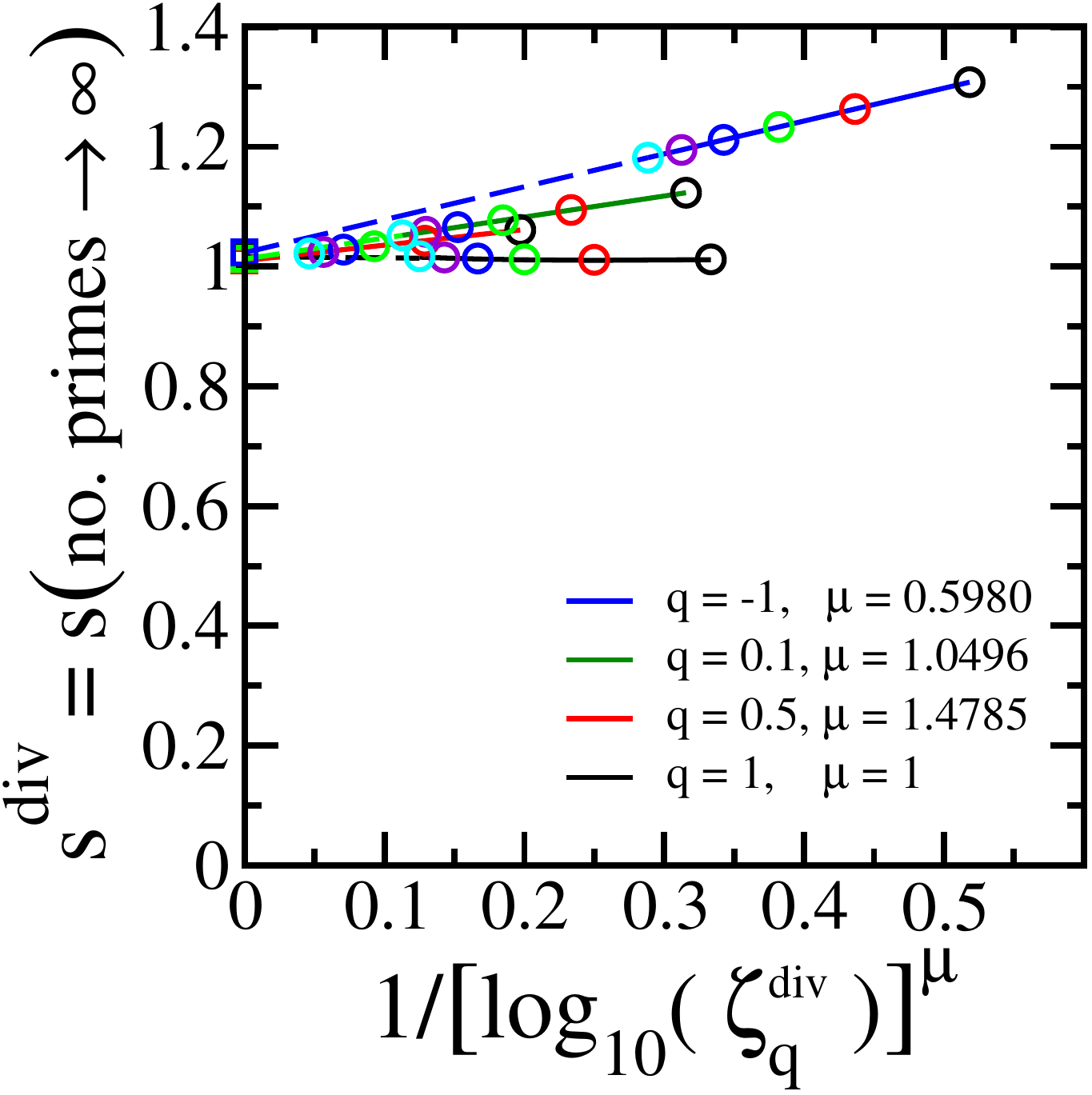}
\caption{\label{fig:coef_lin_1sulog10nprimes.to.sigma.to.mu} 
         \small
         The same as Fig.\ \ref{fig:zeta1}(top right and bottom right)
         for typical values of $q \le 1$.
         The colors of the open circles identify the value of 
         $\zeta_q^{\text{div}}$
         as used in Fig.\ \ref{fig:zeta1}(top middle and bottom middle).
         The colors of the solid lines identify the value of $q$
         according to the legend in the present figure.
         The divergence of $\zeta_q(s)$ occurs at $s^{\text{div}(\infty)}=1$
         (open squares)
         within a numerical error less than 3\%
         for the values of $q$ that we have checked.
        }
\end{figure}

Fig.\ \ref{fig:coef_lin_1sulog10nprimes.to.sigma.to.mu}
shows the final step of the procedure
(see Fig.\ \ref{fig:zeta1} top right and bottom right panels)
for different values of $q \le 1$,
and the maximal estimated numerical error 
is less than 3\% for the values of $q$ that we have checked. 
This result definitely differs from what a brief glance at 
Fig.\ \ref{fig:q-zeta1} might induce one to think.

\begin{figure}
    \centering
    \includegraphics[width=0.30\textwidth]
     {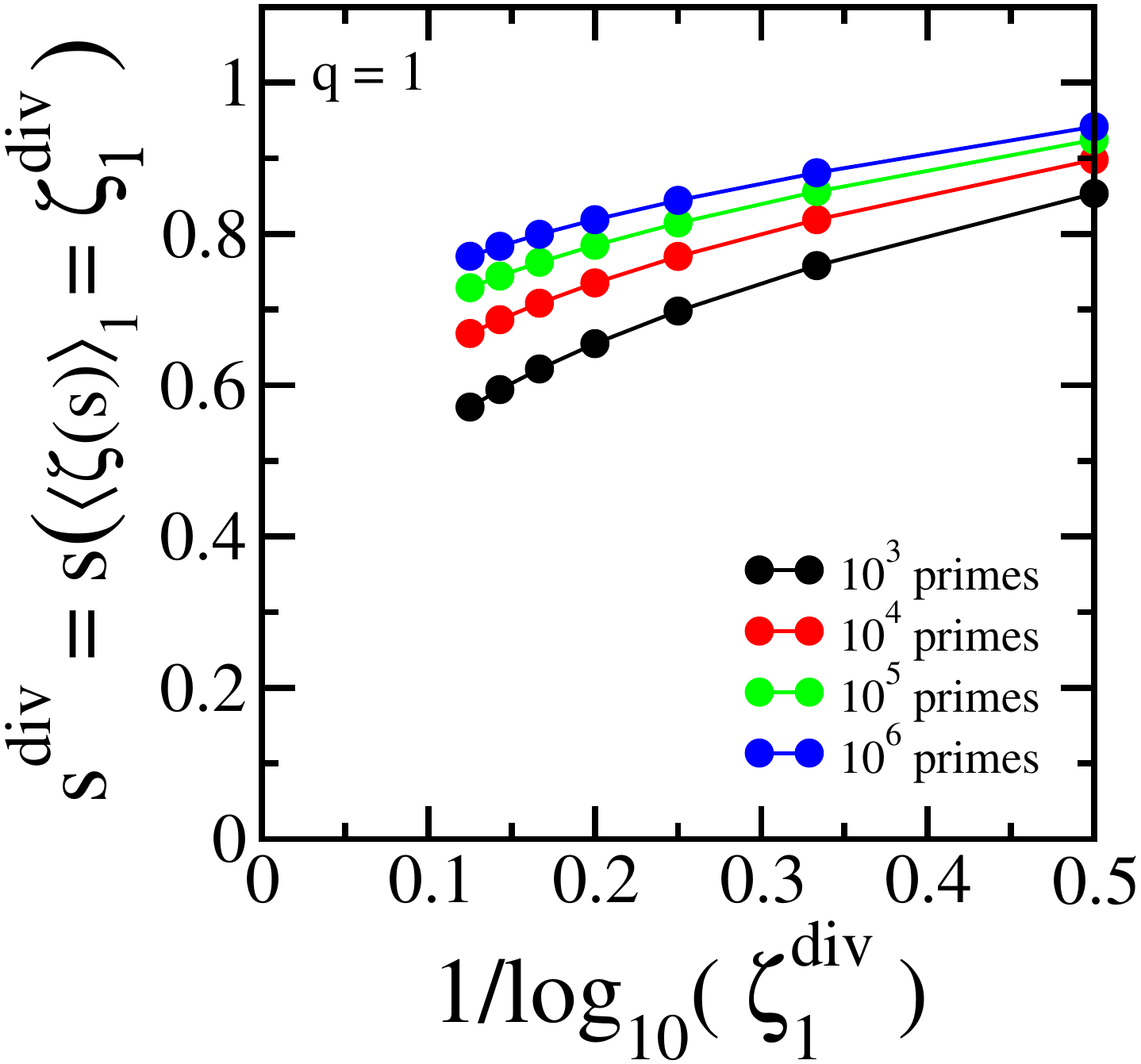}
    \includegraphics[width=0.30\textwidth]
     {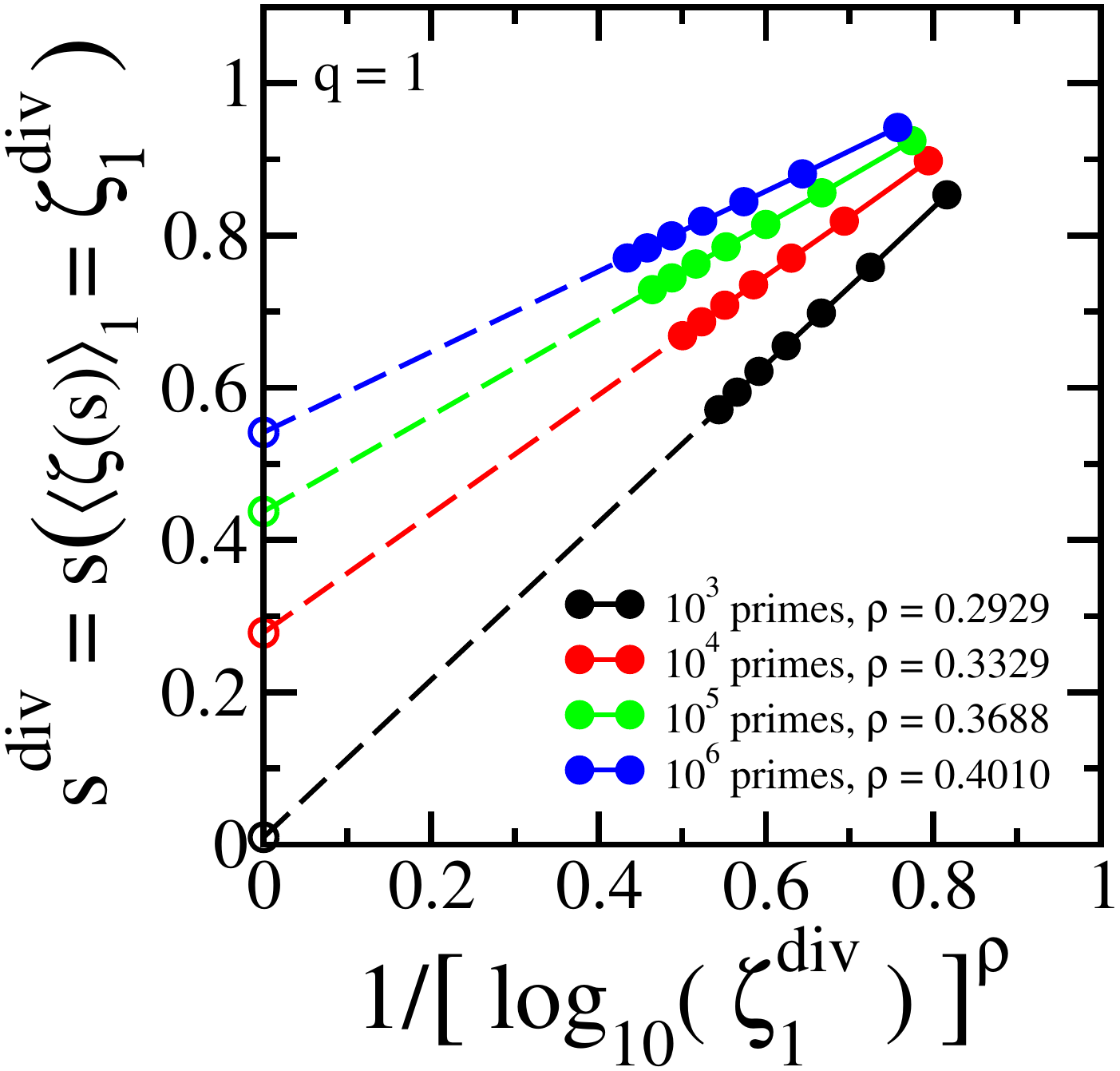}
    \includegraphics[width=0.32\textwidth]
    {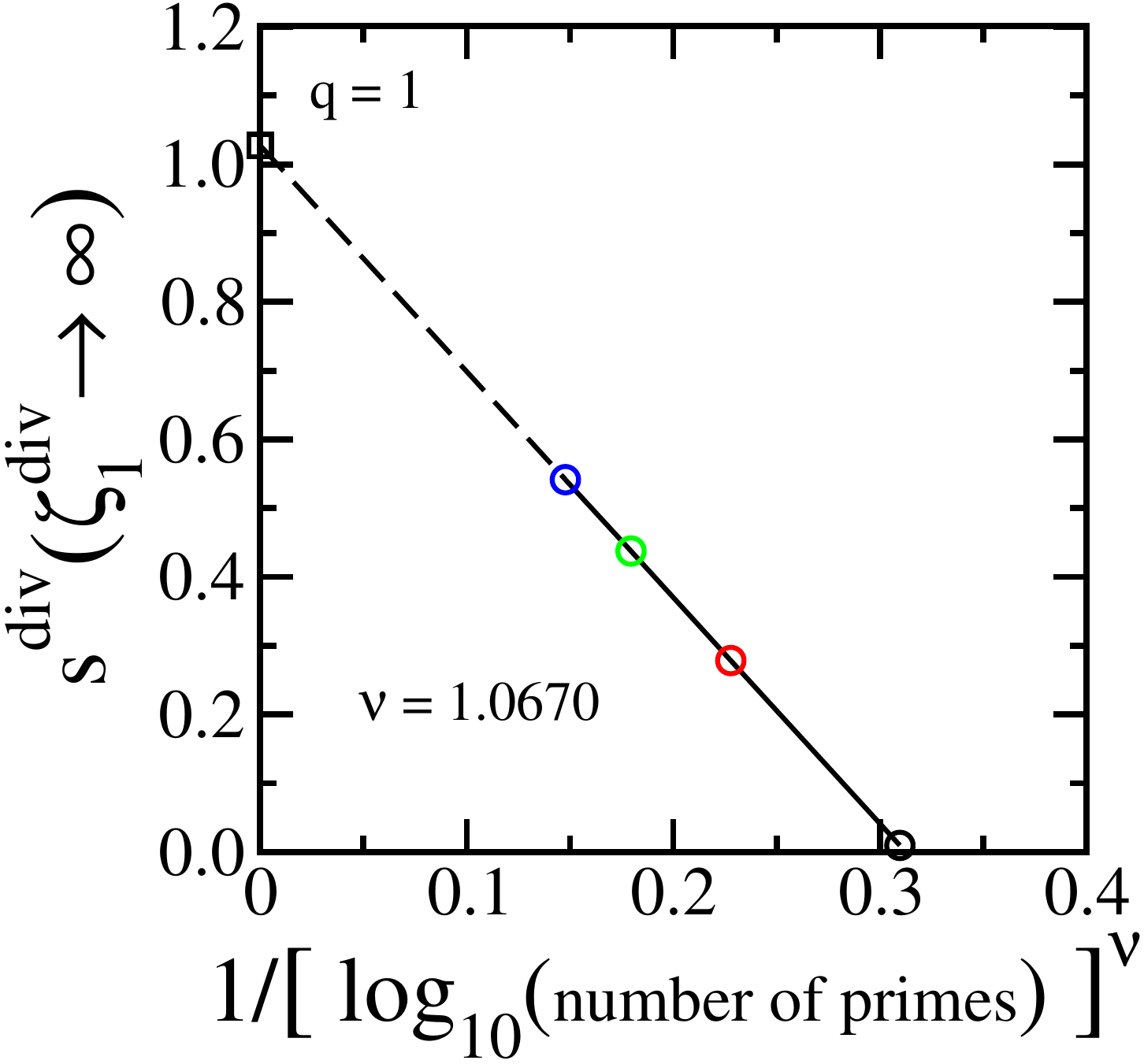}

    \includegraphics[width=0.30\textwidth]
     {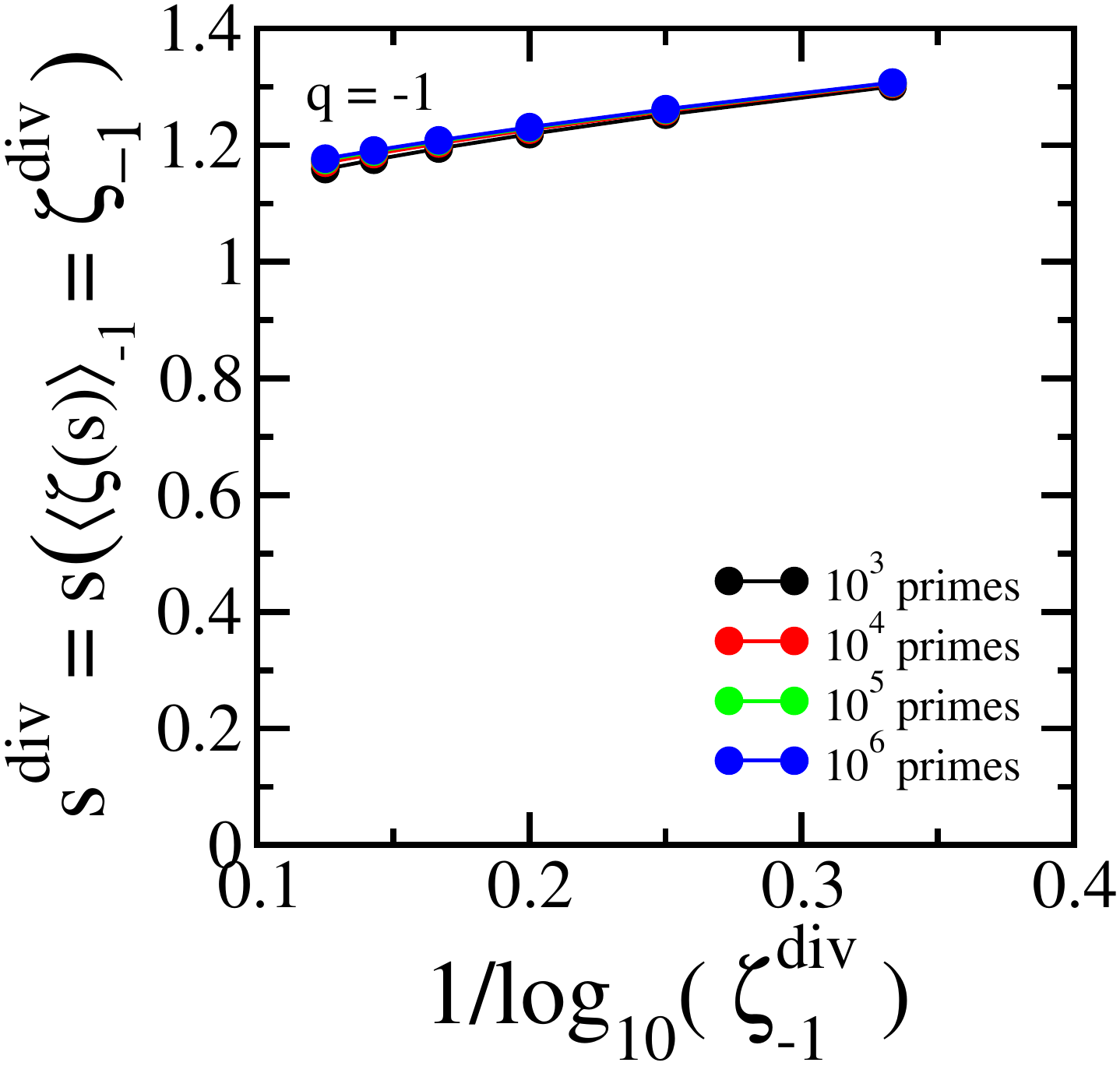}
    \includegraphics[width=0.30\textwidth]
     {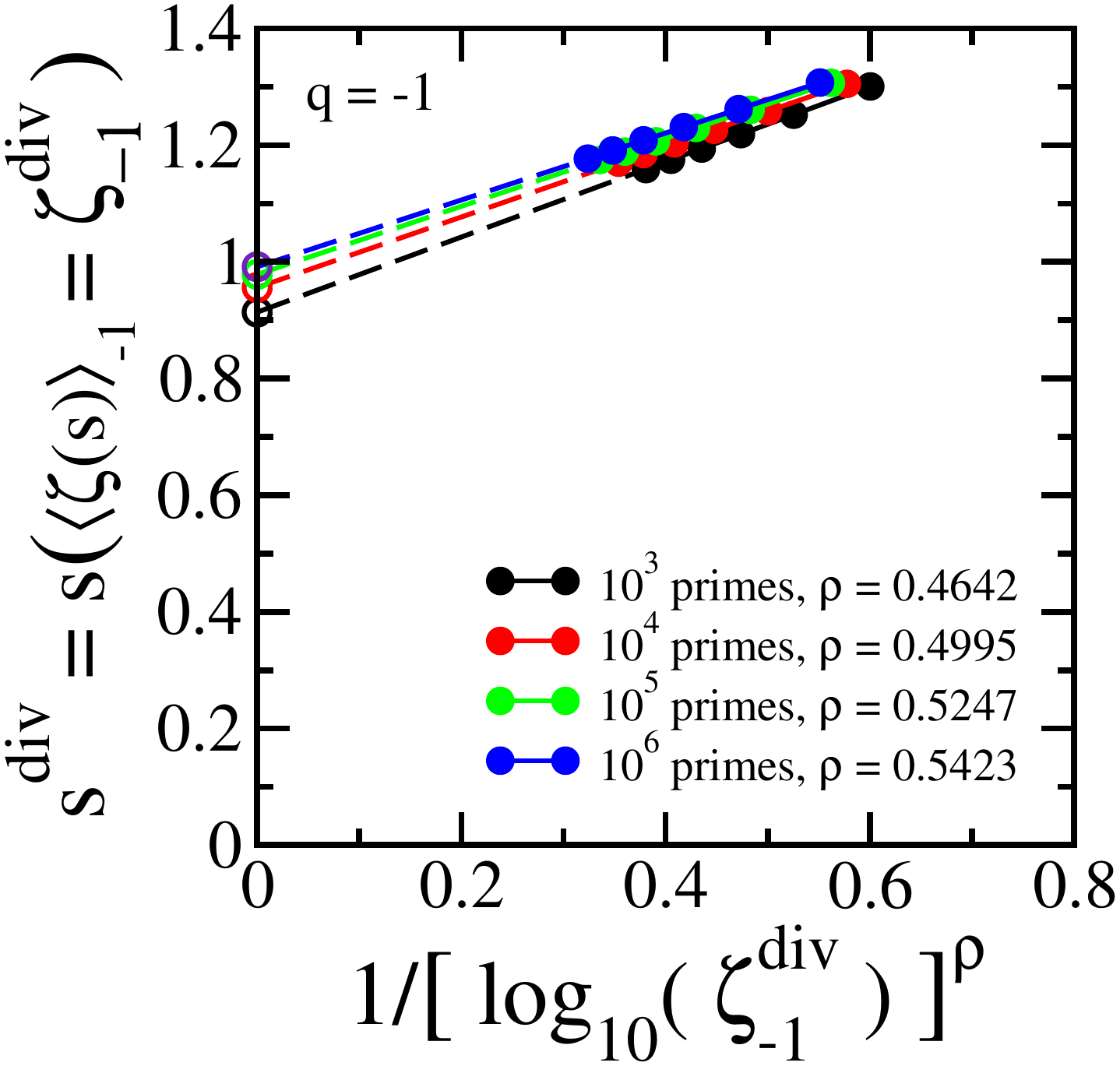}
    \includegraphics[width=0.32\textwidth]
    {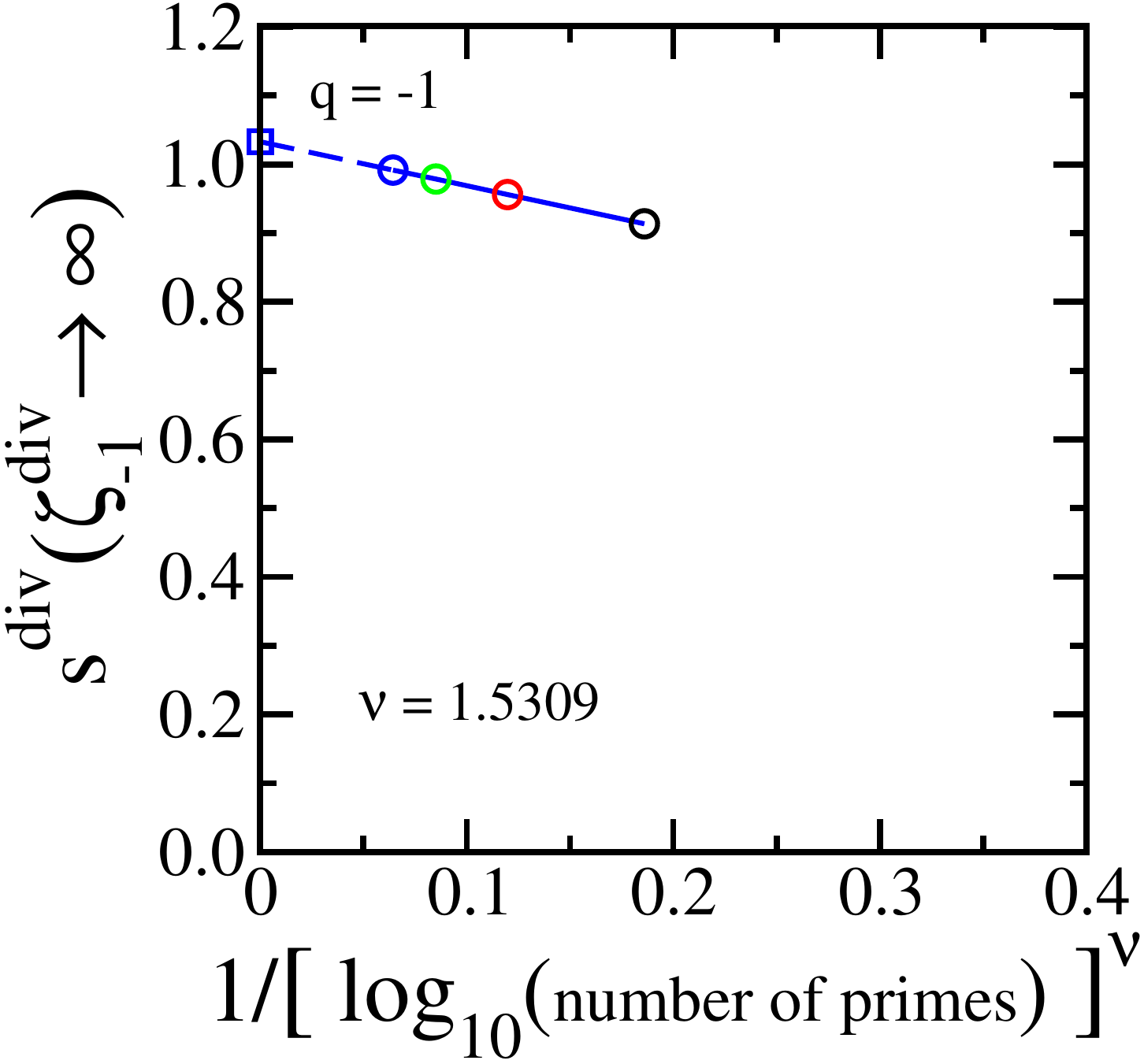}
\caption{\label{fig:zeta1-bis} 
         \small
         {\it Top left panel}:
         Dependence of  
         $s^{\text{div}}=s(\zeta_1^{\text{div}})$ with 
         $1/\log_{10}(\zeta_1^{\text{div}})$
         in Eq.\ (\protect\ref{zeta-2});
         {\it Top middle panel}:
         Power-law rescaling of 
         $s^{\text{div}}=s(\zeta_1^{\text{div}})$ with
         $1/\left[\log_{10}(\zeta_1^{\text{div}})\right]^\rho$.
         The curves point towards $s^{\text{div}}$ 
         with $\zeta_q^{\text{div}}\to\infty$ (open circles);
         {\it Top right panel}:
         Extrapolated values of Fig.\ \ref{fig:zeta1-bis}(top middle) 
         linearly rescaled with 
         $1/\left[\log_{10}(\text{number of primes})\right]^\nu$
         point towards
         $s^{\text{div}(\infty)}=1$ (open square)
         ($ \lim_{\text{number of primes}\to\infty}
            \lim_{\zeta_1^{\text{div}}\to\infty} s^{\text{div}} = 1$)
         within a numerical error less than 3\%.
         The colors of the open circles refer to the number of primes
         identified in Fig.\ \ref{fig:zeta1-bis}(top middle).
         {\it Bottom panels}: the same as top panels with $q=-1$,
         Eq.\ (\protect\ref{eq:zeta_q.eq.q-zeta_1}). 
        }
\end{figure}

\begin{figure}
    \centering
    \includegraphics[width=0.30\textwidth]
     {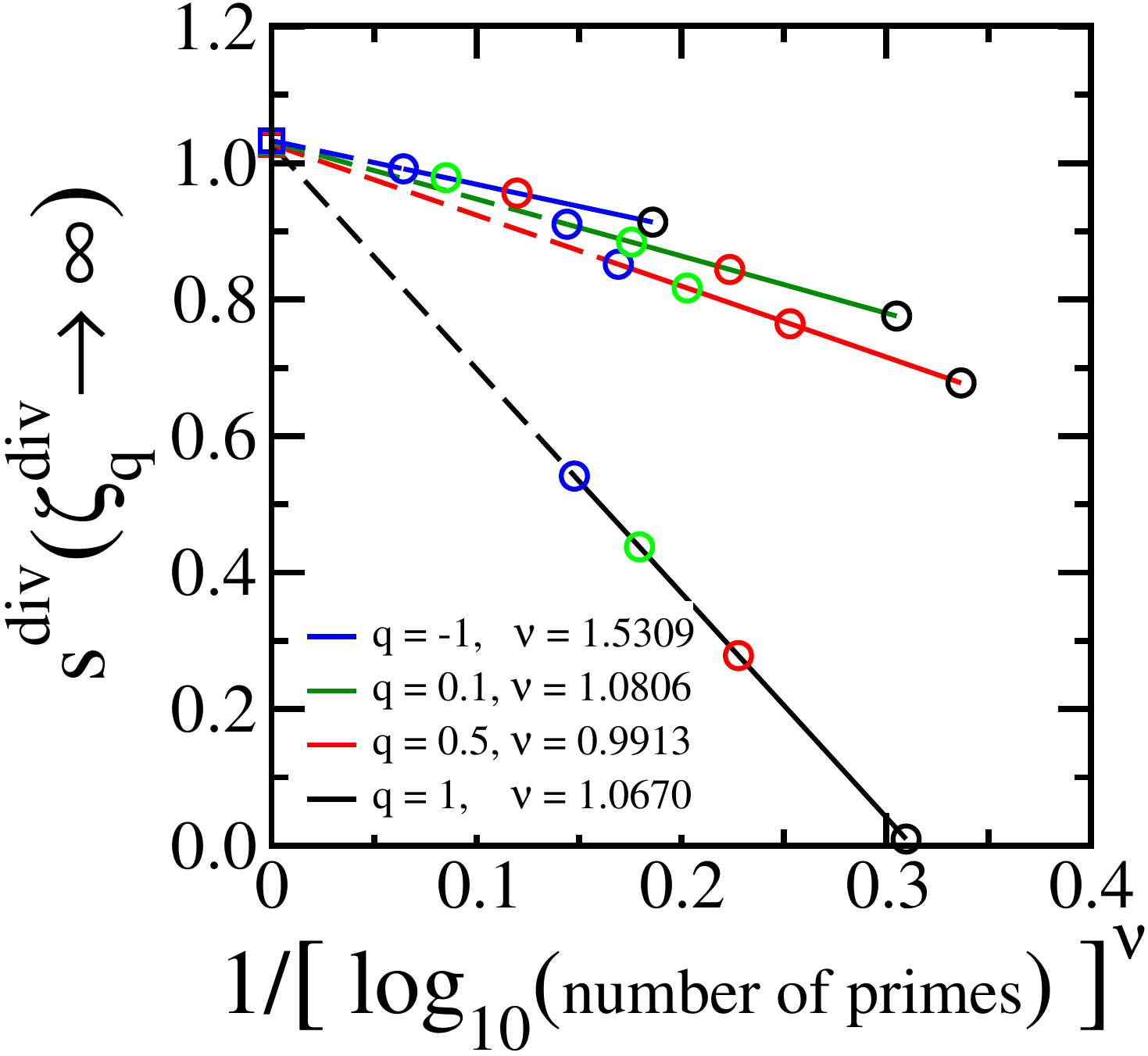}
\caption{\label{fig:coef_lin_1sulog10nprimes.to.rho.to.nu}
         \small
         The same as Fig.\ \ref{fig:zeta1-bis}(top right and bottom right) 
         for typical values of $q \le 1$.
         The colors of the open circles identify the number of primes
         as used in Fig.\ \ref{fig:zeta1-bis}(top middle and bottom middle).
         The colors of the solid lines identify the value of $q$
         according to the legend in the present figure.
         The divergence of $\zeta_q(s)$ occurs at $s^{\text{div}(\infty)}=1$
         (open squares)
         within a numerical error less than 4\%
         for the values of $q$ that we have checked.
        }
\end{figure}

The numerical procedure can be taken in the inverse order,
taking $\zeta_q^{\text{div}}\to\infty$ as the first step,
and then taking increasing number of primes:
see Figs.\ \ref{fig:zeta1-bis} 
and \ref{fig:coef_lin_1sulog10nprimes.to.rho.to.nu}.
Each curve in Fig.\ \ref{fig:zeta1-bis}(top left and bottom left)
displays the value of $s^{\text{div}}$ calculated 
with the same number of primes in Eq.\ (\ref{eq:zeta_q.eq.q-zeta_1}) 
($10^3, 10^4, 10^5, 10^6$ primes)
as a function of $1/\log_{10} \zeta_q^{\text{div}}$.
Here, similarly to Fig.\ \ref{fig:zeta1}(top left and bottom left), 
the curves are not straight lines, so they can hardly be extrapolated.
The empirical power-law rescaling shown in 
Fig.\ \ref{fig:zeta1-bis}(top middle and bottom middle)
indicates that $s^{\text{div}}$ linearly scales with 
$1/[\log_{10} \zeta_q^{\text{div}}]^\rho$
($\rho$ depends on the number of primes),
and the generated straight lines point to the corresponding values of 
$s^{\text{div}}$ with infinite number of primes in the $\zeta_q(s)$ function
(open circles of the top middle and bottom middle panels).
These extrapolated values are rescaled according to a power-law
shown in Fig.\ \ref{fig:zeta1-bis}(top right and bottom right), 
with the empirical power $\nu$ depending on $q$.

Fig.\ \ref{fig:coef_lin_1sulog10nprimes.to.rho.to.nu}
is equivalent to top right and bottom right panels of Fig.\ \ref{fig:zeta1-bis}
for different values of $q \le 1$.
All these cases indicate
$ \lim_{\text{number of primes}\to\infty} 
  \lim_{\zeta_q^{\text{div}}\to\infty} 
  s^{\text{div}} = 1$
within a numerical error less than 4\% 

The empirical powers ($\sigma$, $\rho$, $\mu$, $\nu$) 
have been estimated by fitting a parabola $y = a + bx + c x^2$ 
to the corresponding curve, 
$x$ is the variable of the abscissa of the corresponding 
middle and right panels,
$y = s^{\text{div}}$,
and the fitting value of the power is that for which $c \approx 0$,
estimated with four digits for the power parameter.
The coefficient $a$ of the fitting of the parabola 
is the extrapolated value of the corresponding curve
(open circles of the top middle and bottom middle panels
 of Figs.\ \ref{fig:zeta1} and \ref{fig:zeta1-bis},
 open squares of the top right and bottom right panels
 of Figs.\ \ref{fig:zeta1} and \ref{fig:zeta1-bis}
 and open squares of Figs.\ 
 \ref{fig:coef_lin_1sulog10nprimes.to.sigma.to.mu} 
 and
 \ref{fig:coef_lin_1sulog10nprimes.to.rho.to.nu}).

Similar behavior
is expected for $\zeta_q(s)$ evaluated with the version 
with summations, Eq.\ (\ref{eq:zeta-q-sums}).


\section{\label{non-factorizing-q-primes}
         Algebra violating factorizability of $q$-integer numbers
         in $q$-prime numbers}

The $q$-product corresponding to Eq.\ (\ref{eq:operators-b}),
\begin{equation}
 \langle x \otimes_q  y \rangle_q = \langle x \rangle_q \, \langle y \rangle_q,
\end{equation}
is defined as%
\footnote{Eq.\  (7) of \protect\cite{Borges2004},
          Eq.\ (48) of \protect\cite{BorgesCosta2021}.
         }
\begin{equation}
 x \, \otimes_q \, y  \equiv \Bigl[x^{1-q}+y^{1-q}-1\Bigr]^{\frac{1}{1-q}} 
 \;\;\;\;
 (x \ge 0, \,  \, y \ge 0;\, x \otimes_1  y=xy) 
\label{qproduct}
\end{equation}
or, equivalently,
\begin{equation}
 \label{eq:q-product-b}
 x \otimes_q  y  \equiv e_q^{\ln_q x + \ln_q y} \,,
\end{equation} 
and the $q$-sum is defined as%
\footnote{Eq.\ (A19) of \protect\cite{BorgesCosta2021}.}
\begin{eqnarray}
 \label{eq:q-sum-b}
x \oplus_{q} y &\equiv& e_q^{ \ln \, [e^{\ln_q x} + e^{\ln_q y} ]}  
\nonumber \\
                    &=& \Bigl\{1+(1-q)\ln \Bigl[
                                                  e^{\frac{x^{1-q}-1}{1-q}}
                                                + e^{\frac{y^{1-q}-1}{1-q}}
                                           \Bigr] 
                         \Bigr\}^{1/(1-q)}.
\end{eqnarray}
If $q\ne1$, 
$\langle x \rangle_q \langle y \rangle_q \ne \langle x\,y \rangle_q$.
Eq.\ (\ref{eq:q-sum-b}) with the symbol $\oplus_q$ 
is equivalent to Eq.\ (44) and (A19) of \cite{BorgesCosta2021}
with the symbol ${}_{\{q\}}\oplus$, called oel-addition%
\footnote{
          The notation $\oplus_q$ adopted in Eq.\ (\ref{eq:q-sum-b})
          was used as Eq.\ (4) of \cite{Borges2004}
          with a {\it different} meaning than here,
          namely $x \oplus_qy \equiv x+y + (1-q)xy$,
          which is usually referred to as $q$-sum.
          $x+y + (1-q)xy$ is denoted in Ref.\ \cite{BorgesCosta2021}
          (Eq.\ (25))
          with the symbol 
          $x \, {}_{[q]}\!\oplus y$
          and is called ole-addition.
          }.

The properties (\ref{closure-q_sup})--(\ref{neutral-q_sup-prod})
also hold for the $\Circle_q$ operations (see \cite{BorgesCosta2021}).

\begin{equation}
\zeta^{\Sigma^{\prime}}_q(s) \equiv
 \overset{\!\!\infty}{{\underset{\!\!n=1}{\sideset{}{_q}\sum}}} \; \frac{1}{n^s}
  = 1 \oplus_{q} \frac{1}{2^s} \oplus_{q} \frac{1}{3^s} \oplus_{q} \dots
 \label{zetaprimeprimesum}
\end{equation}
and
\begin{eqnarray}
 \zeta^{\Pi^{\prime}}_q(s) 
  &\equiv&
           \overset{}{{\underset{p\,prime}{\sideset{}{_q}\prod}}} 
           \frac{1}{1-p^{-s}}
\nonumber\\
           &=& \frac{1}{1-2^{-s}} 
               \otimes_q \frac{1}{1-3^{-s}} 
               \otimes_q \frac{1}{1-5^{-s}} 
               \otimes_q \cdots
 \label{zetaprimeprimeproduct}
\end{eqnarray}

In the next Section we present details on a specific generalization.

\section{$q$-Generalizations of the $\zeta(s)$ function
         directly stemming from $q$-numbers}

\begin{figure}[h!]
\includegraphics[width=7.2cm]{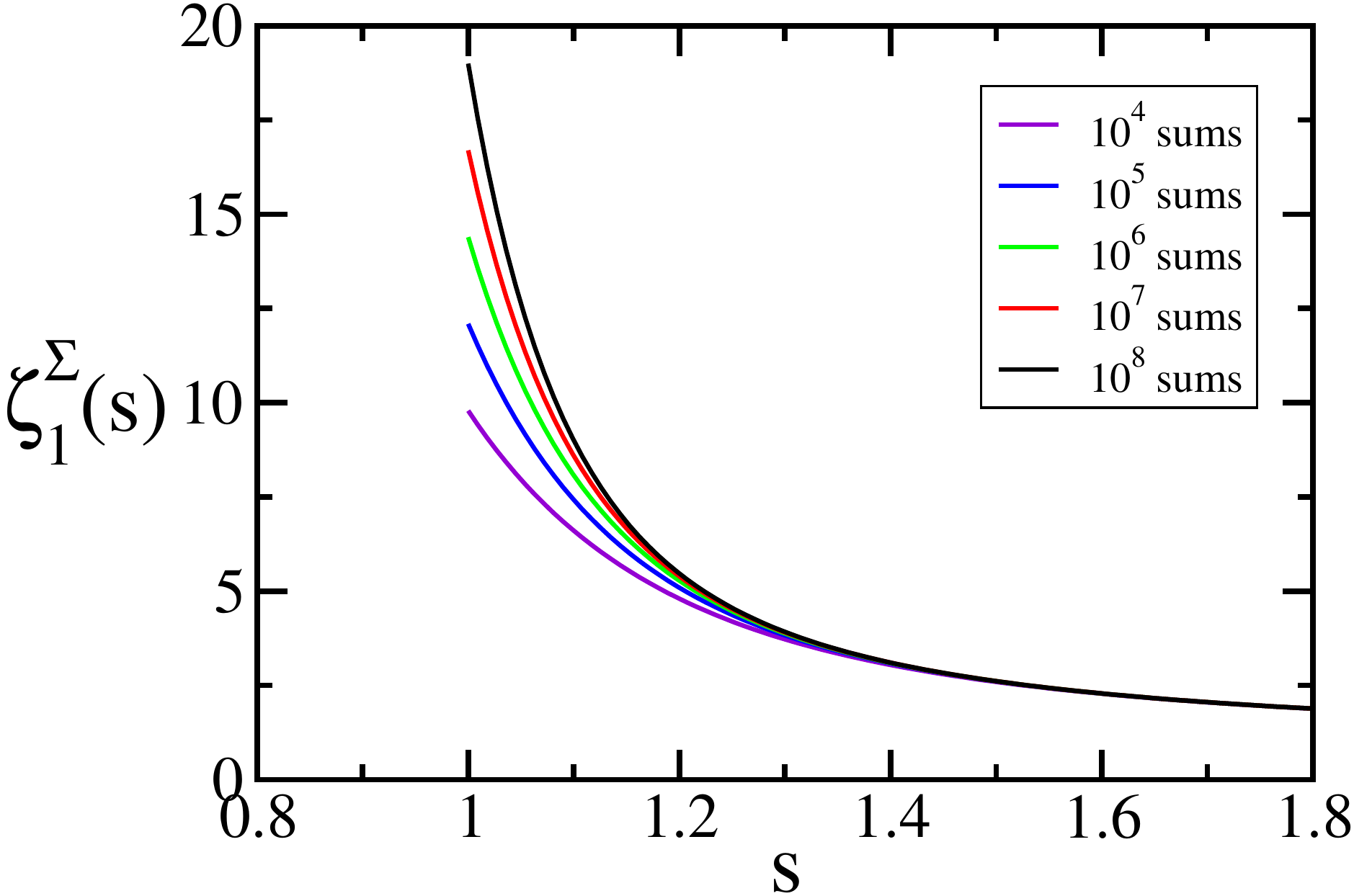}
\includegraphics[width=7.2cm]{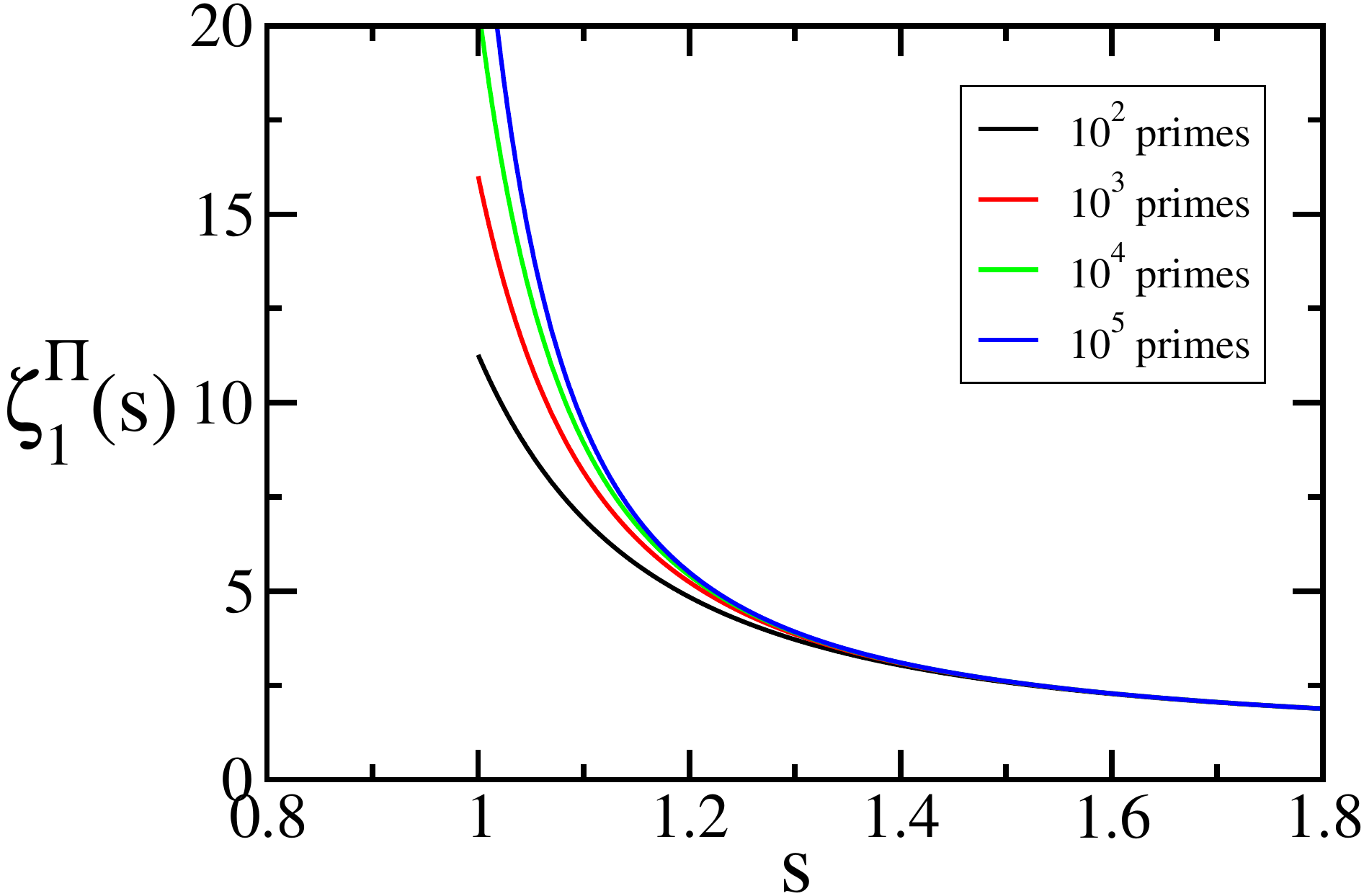}
\centering
\includegraphics[width=7.2cm]{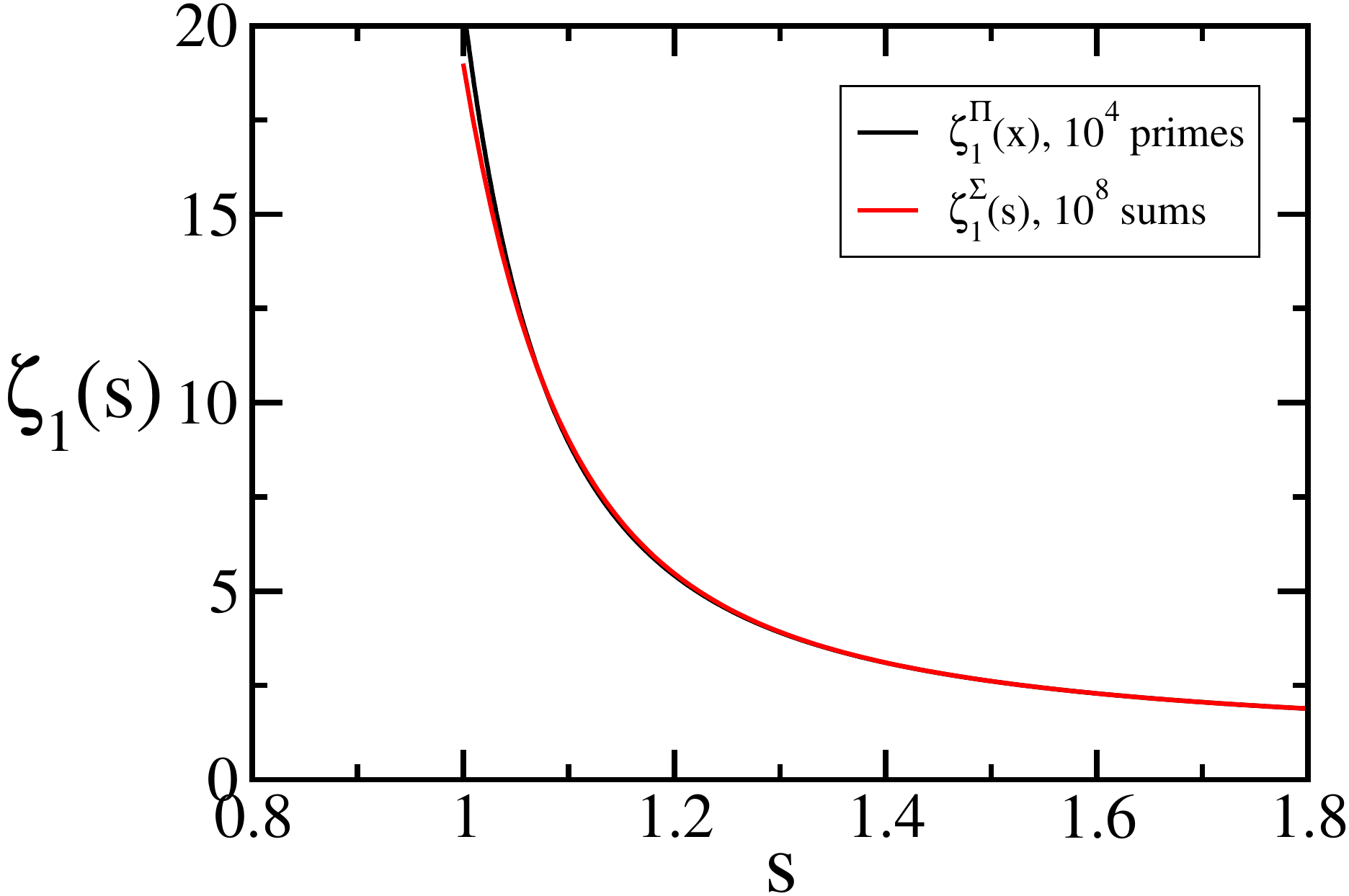}
\caption{\label{fig2}
         \small
         Influence on $\zeta^\Sigma_1(x)=\zeta^\Pi_1(x)=\zeta(x)$ 
         of the number of terms in the sum and in the product, where we have used
         respectively Eq.\ (\ref{zetasum}) and (\ref{zetaproduct}). 
        }
\end{figure}

\begin{figure}
\includegraphics[width=6.6cm]{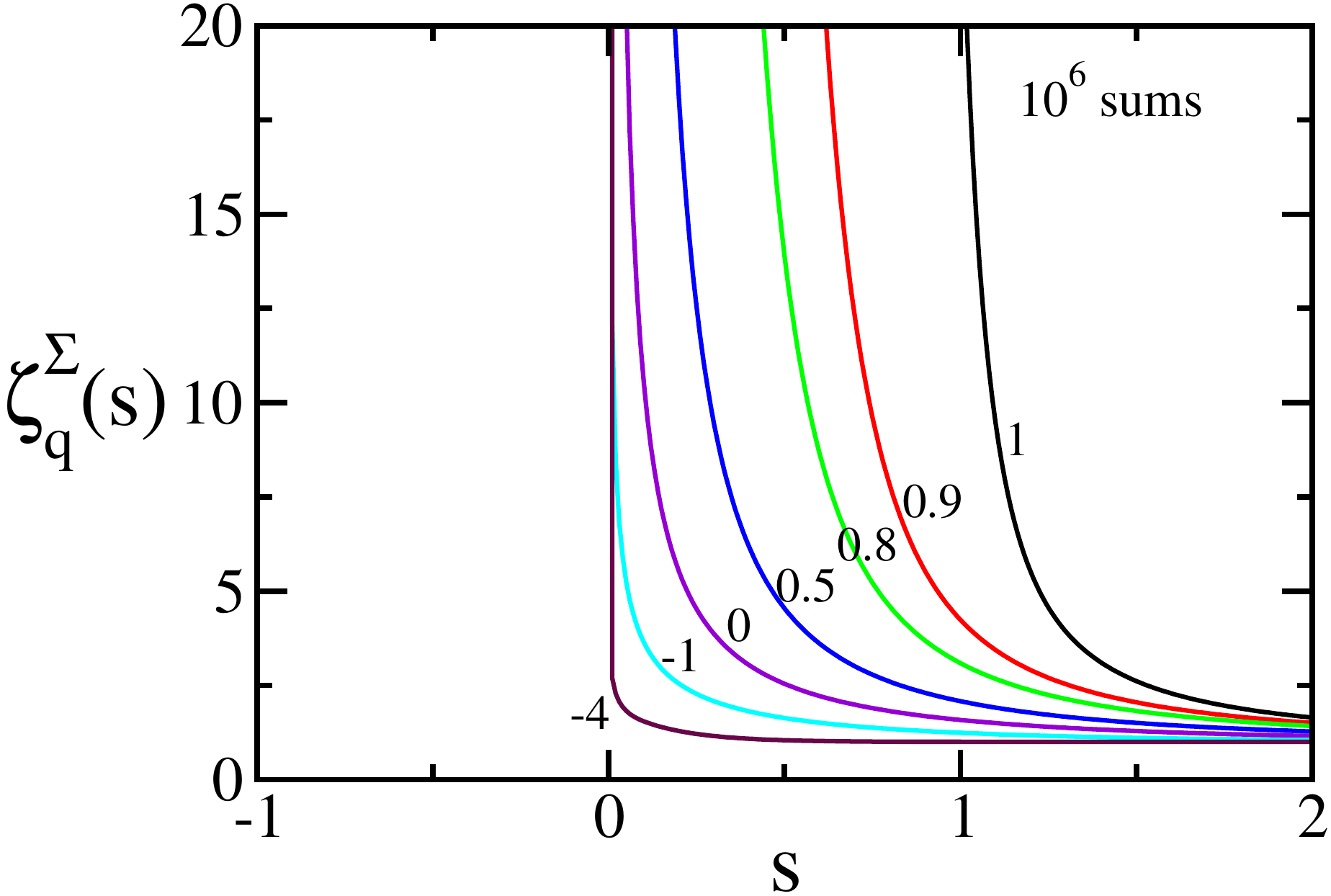}
\includegraphics[width=6.6cm]{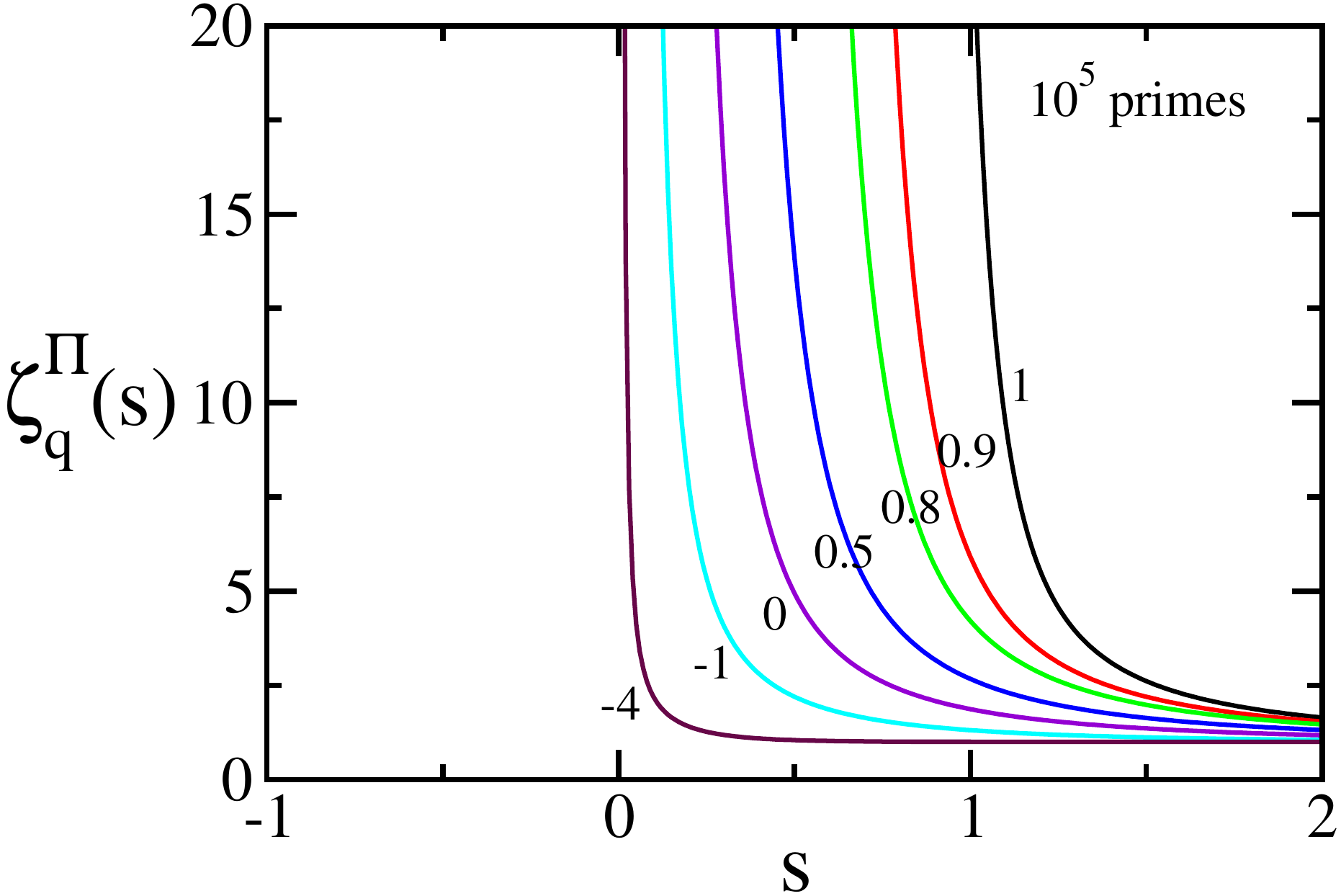}
\centering
\includegraphics[width=6.6cm]{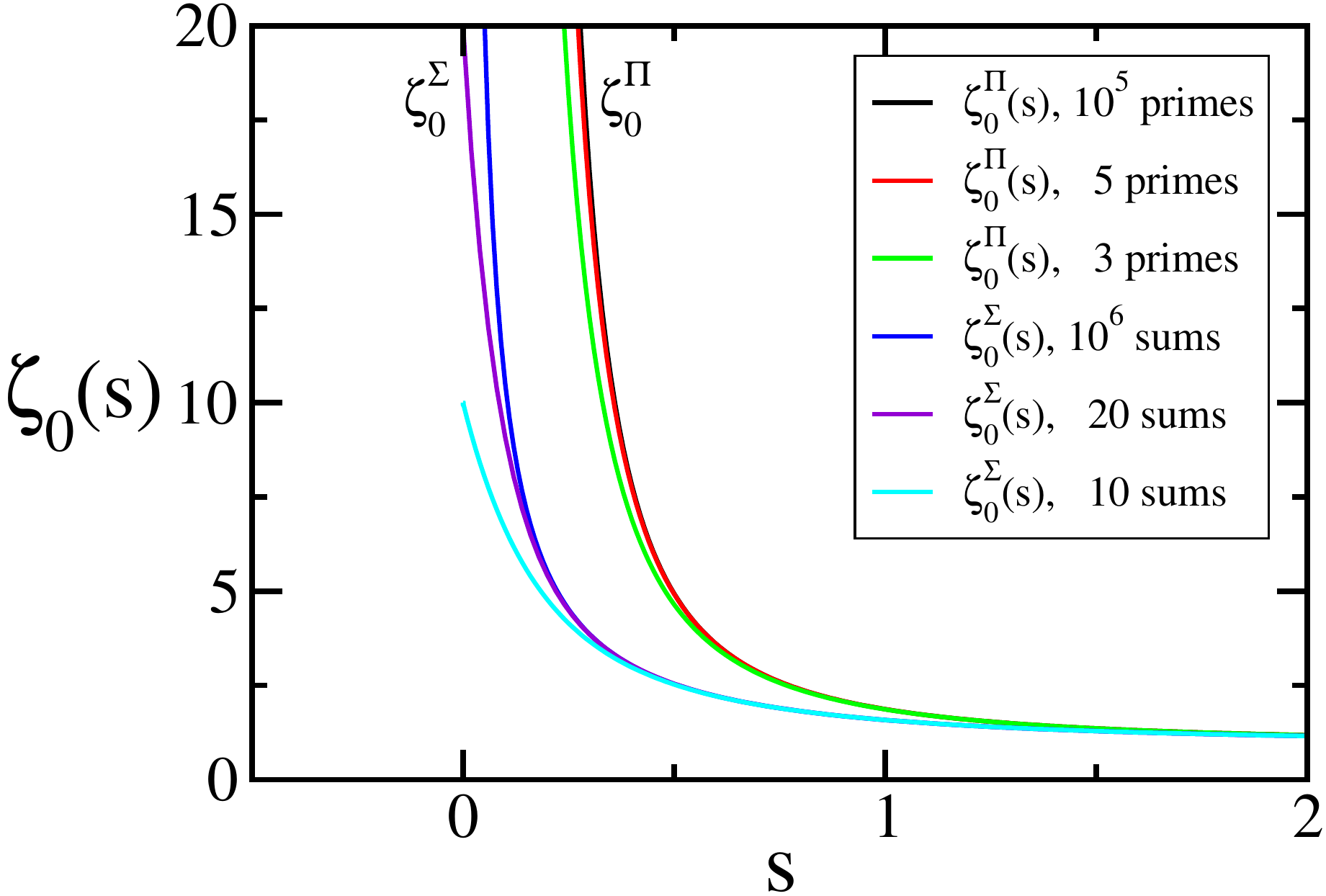}
\caption{\label{fig3} 
         \small
         $\zeta^\Sigma_q(s)$ and $\zeta^\Pi_q(s)$ 
         for typical values of $q \le 1$. 
        }
\end{figure}

We define the following $q$-generalizations of the Riemann $\zeta$ function:
\begin{equation}
 \zeta^\Sigma_q(s) \equiv  \sum_{n=1}^\infty \frac{1}{\langle n \rangle_q^s}
                   =       1 + \frac{1}{\langle 2 \rangle_q^s} +
                               \frac{1}{\langle 3 \rangle_q^s} + \dots 
 \;(s \in \mathbb{R})
 \label{zetasum}
\end{equation}
and 
\begin{eqnarray}
 \zeta^\Pi_q(s) &\equiv& \prod_{p\,prime} \frac{1}{1-\langle p \rangle_q^{-s}}
 \nonumber\\
                &=& \frac{1}{1-\langle 2 \rangle_q^{-s}}
                    \frac{1}{1-\langle 3 \rangle_q^{-s}}
                    \frac{1}{1-\langle 5 \rangle_q^{-s}}
                    \cdots \, ,
 \label{zetaproduct}
\end{eqnarray}
where $\langle n \rangle_q$ is the $q$-number defined in 
Eq.\ (\ref{eq:e_lnq_x}); see Figs.\ \ref{fig2} and \ref{fig3}.

Other possibilities may naturally be considered for generalizing $\zeta(s)$, 
for instance
\begin{equation}
\zeta^{\Sigma^{\prime \prime}}_q(s) 
                \equiv \sum_{n=1}^\infty \Bigl\langle\frac{1}{n^s}\Bigr\rangle_q
                = 1 + \Bigl\langle\frac{1}{2^s}\Bigr\rangle_q 
                    + \Bigl\langle\frac{1}{3^s}\Bigr\rangle_q
                    + \dots 
 \label{zetaprimesum}
\end{equation}
and
\begin{eqnarray}
 \zeta^{\Pi^{\prime \prime}}_q(s) 
            &\equiv& \prod_{p\,prime} \Bigl\langle 
                                        \frac{1}{1-p^{-s}} 
                                       \Bigr\rangle_q
 \nonumber\\
            &=& \Bigl\langle \frac{1}{1-2^{-s}} \Bigr\rangle_q   
                \Bigl\langle \frac{1}{1-3^{-s}} \Bigr\rangle_q
                \Bigl\langle \frac{1}{1-5^{-s}} \Bigr\rangle_q
                \cdots
 \label{zetaprimeproduct}
\end{eqnarray}

Let us however clarify that it is not the scope of the present paper 
to systematically study all such possibilities.

\section{Final remarks}

\begin{figure}[h!] 
 \includegraphics[width=7.9cm]{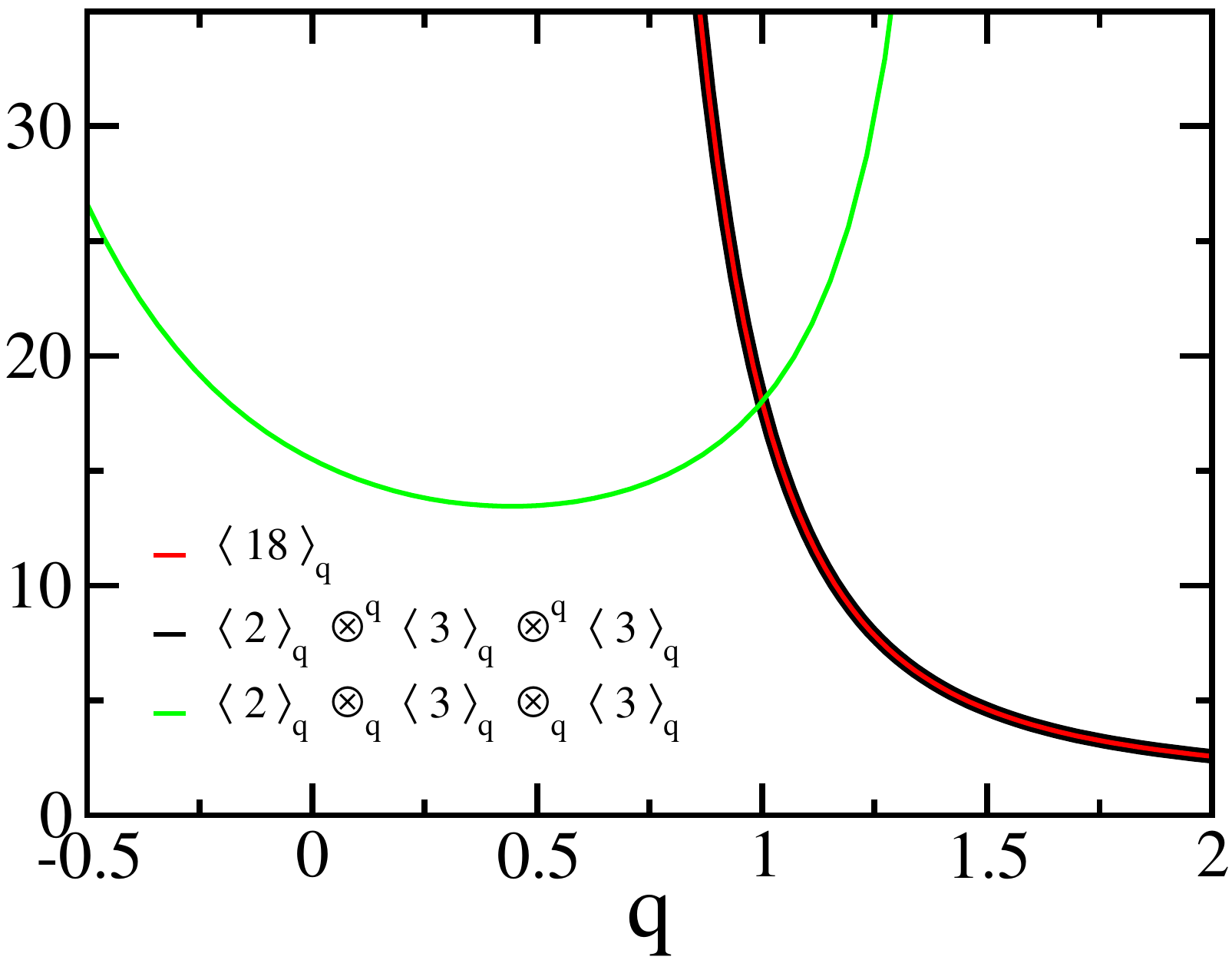}
 \qquad
 \includegraphics[width=7.9cm]{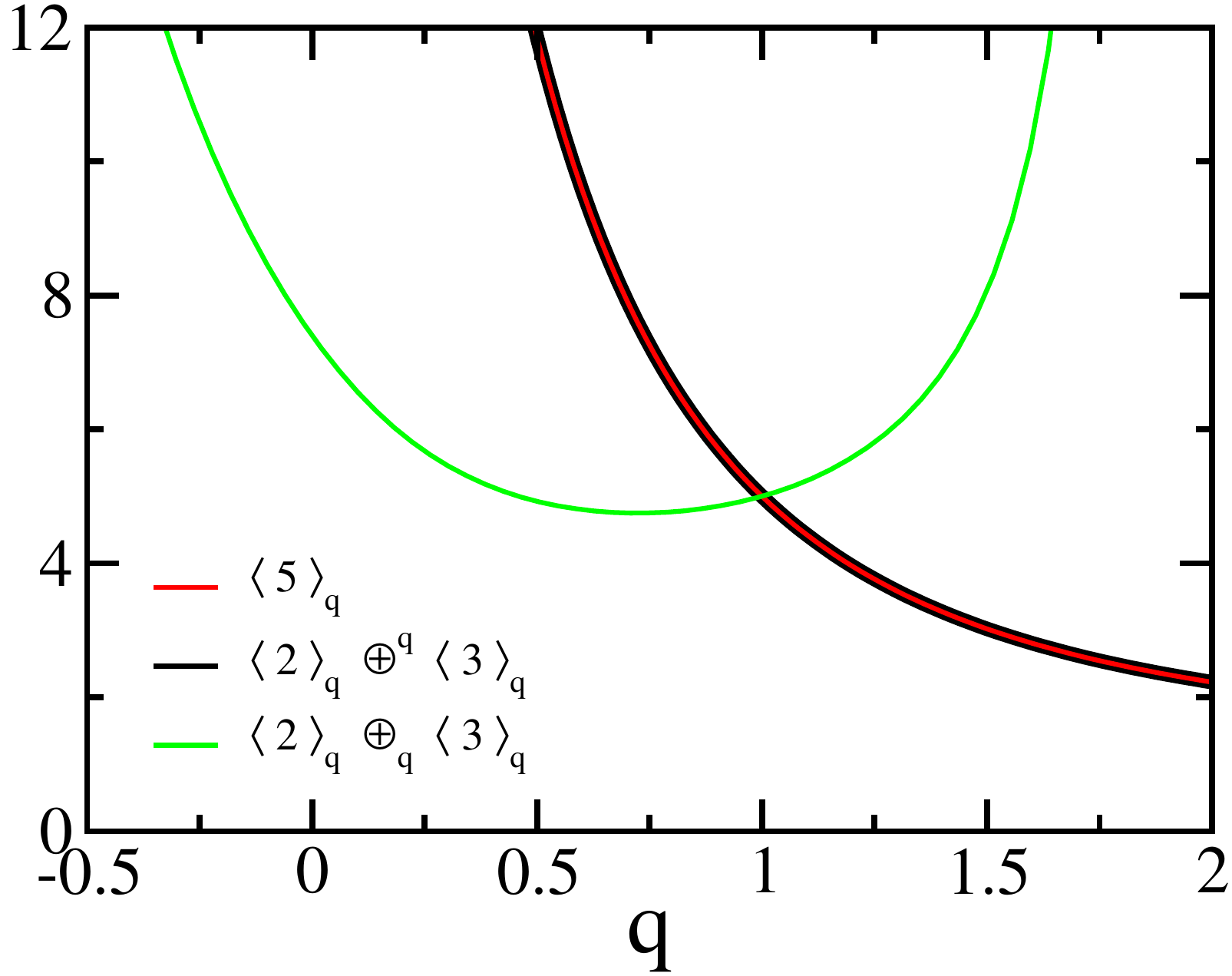}
\caption{\label{18-fact} 
         \small
         {\it Left panel}: 
         Factorization is preserved through the generalized product defined by 
         (\ref{eq:q-product-a}), illustrated with 
         $\langle 2 \rangle_q \otimes^q \langle 3 \rangle_q \otimes^q \langle 3 \rangle_q$
         for different values of $q$ (black curve); 
         the generalized product defined by (\ref{eq:q-product-b}), 
         $\langle 2 \rangle_q \otimes_q \langle 3 \rangle_q \otimes_q \langle 3 \rangle_q$,
         does not preserve factorization (green curve). 
         The red curve is the corresponding $q$-number, Eq.\ (\ref{eq:e_lnq_x}).
         {\it Right panel}: instance of the $q$-sum of $q$-numbers, 
         with the $q$-sums given by (\ref{eq:q-sum-a}),
         $\langle 2 \rangle_q \oplus^q \langle 3 \rangle_q $, (black curve) 
         and (\ref{eq:q-sum-b}), 
         $\langle 2 \rangle_q \oplus_q \langle 3 \rangle_q $, (green curve). 
         The red curve is the corresponding $q$-number, Eq.\ (\ref{eq:e_lnq_x}).
         Notice that there exists a nontrivial value of $q \ne 1$ for which
         $\langle 2 \rangle_q 
          \otimes_q 
          \langle 3 \rangle_q 
          \otimes_q 
          \langle 3 \rangle_q = 18 $
         and, similarly, 
         $\langle 2 \rangle_q \oplus_q \langle 3 \rangle_q = 5 $.
        }
\end{figure}

Let us now illustrate the two algebraic approaches focused on 
in the present paper 
(see Fig.\ \ref{18-fact}):
\begin{equation}
\langle 18\rangle_q = \langle 2 \rangle_q  \otimes^q
                      \langle 3 \rangle_q \otimes^q
                      \langle 3 \rangle_q
                    = \langle 2\rangle_q \otimes^q
                      \bigl(\langle 3 \rangle_q \owedge^q 2 \bigr)
                    \;\;\;(\forall q)
                    \label{factorization-preserving}
\end{equation}
whereas
\begin{equation}
\langle 18\rangle_q \ne \langle 2 \rangle_q  \otimes_q
                        \langle 3 \rangle_q \otimes_q
                        \langle 3 \rangle_q
                     \ne  \langle 2\rangle_q \otimes_q
                        \bigl( \langle 3 \rangle_q \owedge_q 2 \bigl)
                        \;\;\;(\forall q \ne 1).
                        \label{factorization-violating}
\end{equation}
Analogously we have
\begin{equation}
 \langle 5 \rangle_q = \langle 2 \rangle_q \oplus^q \langle 3 \rangle_q
\;\;\;(\forall q)
\label{addition-preserving}
\end{equation}
whereas
\begin{equation}
 \langle 5 \rangle_q \ne \langle 2 \rangle_q \oplus_q \langle 3 \rangle_q
 \;\;\;(\forall q \ne 1).
 \label{addition-violation}
\end{equation}

The algebra preserving the factorizability of $q$-integer numbers 
into $q$-prime numbers (see equality (\ref{factorization-preserving})) 
achieves this remarkable property essentially because it is isomorphic 
to the usual prime numbers. 
On the other hand, precisely because of that, it is unable to properly 
$q$-generalize the concept of a vectorial space in terms of nonlinearity. 
In contrast, the algebra which violates the factorizability of 
$q$-integer numbers into $q$-prime numbers 
(see inequality (\ref{factorization-violating})), 
or some similar algebra, emerges as a possible path for achieving 
the concept of nonlinear vector spaces, 
which has the potential of uncountable applications in theoretical chemistry 
and elsewhere. 

Since the inequality relation between $q$-primes remains the same 
as that for $q=1$, it is plausible that the nontrivial zeros 
in the analytic extension behaves similarly.
More precisely, it might well be that, by extending 
$\zeta_q(s)$, $\zeta^\Sigma_q(s)$ and $\zeta^\Pi_q(s)$ 
to the complex plane $z$, all nontrivial zeros belong to 
specific single continuous curves, 
$\mathbb{R}(z) = f_q(\mathbb{I}(z))$,
thus $q$-generalizing the $q=1$ Riemann's 1859 celebrated conjecture 
$\mathbb{R}(z)=1/2$.

We have here explored generalizations of the $\zeta(s)$ function
based on a specific type of $q$-number, Eq.\ (\ref{eq:e_lnq_x}), 
and two associated generalized algebras 
(Sections \ref{factorizing-q-primes} and \ref{non-factorizing-q-primes}).
Three additional forms of $q$-generalized numbers,
Eqs.\ (\ref{eq:eq_ln_x})--(\ref{eq:lnq_e_x}), 
are identified in Ref.\ \cite{BorgesCosta2021}.
To each of these $q$-numbers, we can associate two consistently generalized
algebras, one of them violating the factorizability in prime numbers
(see Ref.\ \cite{BorgesCosta2021}), the other one following along the
lines of Section \ref{factorizing-q-primes}.
Similarly, various other generalizations
of the $\zeta(s)$ function may of course be developed.
Naturally, the extension of the present $q$-generalized $\zeta(s)$ functions
to complex $z$ surely is interesting, 
but does not belong to the aim of the present effort.
In any case, the intriguing fact that various infinities 
appear to linearly scale with negative powers of logarithms 
might indicate some general
tendencies.

It is well known that both random matrices  and quantum chaos 
(classically corresponding to strong chaos, 
 i.e., {\it positive} maximal Lyapunov exponent) 
\cite{Berry1986,BerryKeating1999,KeatingSnaith2000,FirkMiller2009,strongchaos} 
are related to the Riemann $\zeta$-function and prime numbers.
On the other hand, both random matrices and strong chaos have been conveniently
$q$-generalized, in \cite{ToscanoVallejosTsallis2004} and  
\cite{WeinsteinLloydTsallis2002,WeinsteinLloydTsallis2004,QueirosTsallis2005} 
respectively. These facts open the door for possible applications 
of the present $q$-generalizations of prime numbers and of the 
$\zeta$-function to $q$-random matrices and to weak chaos 
(classically corresponding to {\it vanishing} maximal Lyapunov exponent, 
 which recovers strong chaos in the $q\to 1$ limit). 
Moreover, connections of the present developments within the realm of the 
theory of numbers, or, more specifically, the theory of prime numbers,
remain, at this stage, out of our scope.
Further work along these lines would naturally be very welcome.

\section*{Acknowledgments}
This work has received partial financial support by CNPq, CAPES and FAPERJ 
(Brazilian agencies).


\end{document}